\documentclass[twoside,12pt]{amsart}
\usepackage{amsfonts}
\usepackage{amsmath}
\usepackage{amssymb}
\usepackage[cp850]{inputenc}
\usepackage[bookmarksnumbered,plainpages,hypertex]{hyperref}
\usepackage{color}

\input{xy}
\xyoption{all}
\newtheorem{theorem}{\sc Theorem}[section]
\newtheorem{proposition}[theorem]{\sc Proposition}

\newtheorem{lemma}[theorem]{\sc Lemma}
\newtheorem{corollary}[theorem]{\sc Corollary}
\theoremstyle{definition}
\newtheorem{definition}[theorem]{\sc Definition}

\newtheorem{example}[theorem]{\sc Example}

\theoremstyle{remark}

\def\ot{\otimes}

\newcommand{\diaguno}{\xymatrix@C=40pt{
  U_T \ar[d]_{g} \ar[r]^{g} & U_TG\ar[d]^{U_T   \Delta^G} \\
  U_TG \ar[r]^{g G} & U_TG G  }\hspace {1.5cm}
  \xymatrix{
  U_T \ar[rr]^{g} \ar[dr]_{U_T}
                &  &    U_TG \ar[dl]^{U_T \varepsilon ^G}    \\
                & U_T                 }}

\numberwithin{equation}{section}

\begin{document}

\title{Pre-torsors and Galois comodules over mixed distributive laws}
\author{Gabriella B\"ohm}
\author{Claudia Menini}
\address{Research Institute for Particle and Nuclear Physics, Budapest,
H-1525 Budapest 114, P.O.B.49, Hungary.}
\address{University of Ferrara, Department of Mathematics, Via Machiavelli 35,
Ferrara, I-44100,Italy}
\date{Nov 2008}

\begin{abstract}
We study comodule functors for comonads arising from mixed
distributive laws. Their Galois property is reformulated in terms
of a (so-called) regular arrow in Street's bicategory
of comonads. Between categories possessing equalizers, we
introduce the notion of a regular adjunction.
An equivalence is proven between the category of
pre-torsors over two regular adjunctions $(N_A,R_A)$ and
$(N_B,R_B)$ on one hand, and the category of
regular comonad arrows $(R_A,\xi)$ from some equalizer preserving
comonad ${\mathbb C}$ to $N_BR_B$ on the other. 
This generalizes a known relationship between pre-torsors over
equal commutative rings and Galois objects of coalgebras.
Developing a bi-Galois theory of comonads, we show that a
pre-torsor over regular adjunctions determines also a second
(equalizer preserving) comonad ${\mathbb D}$  and a co-regular
comonad arrow from ${\mathbb D}$ to $N_A R_A$, such that the
comodule categories of ${\mathbb C}$ and ${\mathbb D}$ are
equivalent.
\end{abstract}
\maketitle

\section{Introduction}

In order to generalize Hopf Galois extensions to Galois extensions
by coalgebras, so called {\em entwining structures} were
introduced in \cite{BrMa}. These are examples of
Beck's {\em mixed distributive laws} \cite{Be}, for a monad $(-)
\ot_k {\bf T}$ induced by an algebra ${\bf T}$, and
a comonad $(-)\ot_k {\bf C}$ induced by a coalgebra ${\bf C}$, on the category
of $k$-modules, for a commutative ring $k$.

The basis of the generalization is an observation that a
right comodule algebra ${\bf T}$ of a Hopf algebra ${\bf H}$ is entwined
with the coalgebra underlying ${\bf H}$. Moreover, ${\bf T}$ is
a comodule for the lifted comonad $(-) \ot_{\bf T} ({\bf T} \ot_k {\bf
H})\cong (-)\ot_k {\bf H}$ on the category of ${\bf T}$-modules
(whose comodules are usually called
Hopf modules). Denoting by ${\bf B}$ the endomorphism algebra of
${\bf T}$ as a Hopf module, the ${\bf H}$-Galois
property of the algebra extension ${\bf B} \subseteq {\bf T}$ can
be formulated as a Galois property of the functor $(-)\ot_{\bf B}
{\bf T}:\mathrm{Mod}\textrm{-}{\bf B} \to
\mathrm{Mod}\textrm{-}{\bf T}$. This latter property means that the canonical
comonad morphism 
$$
(-) \ot_{\bf T} ({\bf T \ot_B T}) \to (-) \ot_{\bf T} ({\bf T} \ot_k {\bf H}),
\qquad m \ot_{\bf T} a'\ot_{\bf B} a \mapsto m \ot_{\bf T} a' \varrho(a)
$$
is an isomorphism, where $\varrho:{\bf T} \to {\bf T}\ot_k {\bf H}$ denotes
the coaction.

Although entwining structures were originally introduced to
develop Galois theory for coalgebras in \cite{BrMa}, this method
turned out to have a wider application. Using mixed distributive
laws of a monad and a comonad on the category of modules over an
arbitrary algebra ${\bf R}$, also Galois extensions by
bialgebroids, and more generally by corings, over ${\bf R}$ fit
this scenario. By these motivations, the first aim of this paper
is to study Galois functors for comonads that arise from arbitrary
mixed distributive laws.

A Galois comodule algebra ${\bf T}$ of a $k$-Hopf
algebra ${\bf H}$, for which the Hopf module
endomorphism algebra ${\bf B}$ of ${\bf T}$ is trivial (i.e. it
contains only multiplication by $k$), is called an ${\bf
H}$-Galois object. As observed by Grunspan \cite{Gr}
and Schauenburg \cite{Scha:Qtor}, \cite{Sch}, \cite{Scha:HbiGal},
{\em faithfully flat} Hopf-Galois objects can be described
equivalently, without explicit mention of the coacting Hopf
algebra ${\bf H}$, in terms of {\em torsors}. A torsor means a certain map
${\bf T} \to {\bf T} \ot_k {\bf T} \ot_k {\bf T}$, from which the {\em flat}
Hopf algebra ${\bf H}$ can be reconstructed uniquely up to isomorphism. (For a
study of the case when ${\bf 
T}={\bf H}$ as algebras, not necessarily faithfully flat over $k$,
see also \cite{Skoda}). This observation was generalized to Hopf Galois
extensions of arbitrary algebras ${\bf B}$ (using the notion of a {\em ${\bf
    B}$-torsor}) in \cite{Scha:ba.nc.base}, to Galois extensions by
bialgebroids (using the notion of an {\em ${\bf A}$- ${\bf B}$-torsor}) in
\cite{Hobst:phd} and \cite{BT}, to Galois extensions by corings (using the
notion of a {\em pre-torsor}) in \cite{BT} and to Galois comodules of corings
arising from entwining structures (using the notion of a {\em bimodule herd})
in \cite{BV}.

Note that the definition of a torsor has the symmetry of reversing
the order of the tensor factors in the codomain. This symmetry
leads to the interesting fact \cite{Sch} that a faithfully flat
torsor ${\bf T}$ determines a second {\em flat} Hopf algebra ${\bf
H}'$, such that ${\bf T}$ is a left ${\bf H}'$-comodule algebra
and an ${\bf H}'$-Galois extension of $k$. Moreover, the Hopf algebras ${\bf
  H}$ and ${\bf H}'$, coacting on ${\bf T}$ on the right and on the left,
respectively, are {\em Morita-Takeuchi equivalent}, i.e. they have equivalent
categories of comodules.

Under sufficiently strong assumptions, the construction of two
corings (over respective base algebras ${\bf A}$ and ${\bf B}$)
from an ${\bf A}$-${\bf B}$ pre-torsor or even from a bimodule
herd, has been carried out in \cite{BT} and \cite{BV},
respectively. However, the resulting corings are not known to be
flat, and their comodule categories do not seem to be equivalent
without further, somewhat intricate assumptions, see \cite[Remark
4.7]{BT}. Placing the problem in a more general categorical
context, in this paper we give an explanation of the origin of
this difficulty. Namely, we show that a pre-torsor (or a finitely
generated projective bimodule herd) that is faithfully flat as a
left module for both base algebras ${\bf A}$ and ${\bf B}$,
determines (uniquely up to natural isomorphisms) two
comonads ${\mathbb C}$ and ${\mathbb D}$ on the category of ${\bf
A}$-, and ${\bf B}$-modules, respectively, whose underlying
functors preserve kernels. These two comonads have equivalent
comodule categories. However, even in the situation when one can
associate an ${\bf A}$-coring ${\bf C}$ and a ${\bf B}$-coring
${\bf D}$ to a pre-torsor or a bimodule herd, as in \cite{BT} or
\cite{BV}, it is not guaranteed that the comonads ${\mathbb C}$
and ${\mathbb D}$ are induced by these corings ${\bf
  C}$ and ${\bf D}$.  This holds
exactly if the corings ${\bf C}$ and ${\bf D}$ are {\em flat} left modules over
their respective base algebras, i.e. in the situation discussed in
\cite[Remark 4.7]{BT}.

The paper is organized as follows.

In Section \ref{sec:preli} we recall some results from category
theory that we need as background.

Section \ref{sec:entw} is devoted to a study of Galois functors (in the sense
of \cite[Definition 4.5]{MW}) for a comonad arising from a mixed distributive
law.
The category of Galois functors (with domain category ${\mathcal B}$) for a
comonad arising from a mixed distributive law (of functors on a category
${\mathcal A}$), is described as a suitable subcategory of a newly defined
category $\mathrm{Arr}({\mathcal A},{\mathcal B})$.
The objects of $\mathrm{Arr}({\mathcal A},{\mathcal B})$
consist of two adjunctions
\begin{equation}\label{eq:two_adj}
\xymatrix{
{\mathcal A}\ar@<2pt>[rr]^-{N_A}&&
{\mathcal T}\ar@<2pt>[ll]^-{R_A} \ar@<-2pt>[rr]_-{R_B}&&
{\mathcal B}\ar@<-2pt>[ll]_-{N_B}
},
\end {equation}
together with a comonad
${\mathbb C}$ on ${\mathcal A}$ and a comonad arrow $(R_A,\xi)$ from ${\mathbb
C}$ to the comonad $N_B R_B$.

Generalizing pre-torsors in \cite{BT} (thus in particular generalizing
Grunspan- Schauenburg torsors), in 
Section \ref{sec:pre-tor} we define pre-torsors over two
adjunctions as in \eqref{eq:two_adj}. A pre-torsor is a natural transformation
$R_A N_B \to R_A N_B R_B N_A R_A N_B$, subject to compatibility conditions
with the units and counits of the adjunctions. We consider a full subcategory
of so called {\em regular} pre-torsors that is shown to be equivalent to a
full subcategory of $\mathrm{Arr}({\mathcal A},{\mathcal B})$.
Since in Section \ref{sec:entw} we describe Galois functors
(of comonads arising from mixed distributive laws) via objects in
$\mathrm{Arr}({\mathcal A},{\mathcal B})$, this equivalence relates
in particular such Galois functors to pre-torsors, provided that they obey our
regularity assumptions.
Note that all these regularity assumptions hold for a
(${\bf B}$-)torsor corresponding to
a faithfully flat Galois extension ${\bf B}\subseteq {\bf T}$ by a Hopf algebra
over a field $k$, i.e. when
$$
\xymatrix{
{\mathcal A}= \mathrm{Vec}\textrm{-}k\ar@<2pt>[rrr]^-{N_A=(-)\ot_k {\bf T}}&&&
{\mathcal T}=\mathrm{Mod}\textrm{-}{\bf T}
\ar@<2pt>[lll]^-{R_A=\mathrm{Hom}_{\bf T}({\bf T},-) }
\ar@<-2pt>[rrr]_-{R_B=\mathrm{Hom}_{\bf T}({\bf T},-)}&&&
{\mathcal B}=\mathrm{Mod}\textrm{-}{\bf B}
\ar@<-2pt>[lll]_-{N_B=(-)\ot_{\bf B} {\bf T}}
} .
$$
More generally, these assumptions hold for bimodule herds corresponding to
(left) faithfully flat Galois comodules for entwining structures.

By the equivalences in Section \ref{sec:pre-tor}, we associate in particular
two comonads ${\mathbb C}$ on ${\mathcal A}$ and ${\mathbb D}$ on ${\mathcal
  B}$ to a pre-torsor satisfying our regularity assumptions.
Generalizing the Morita Takeuchi equivalence of two Hopf algebras in a
faithfully flat bi-Galois object, in the final Section \ref{sec:Mor-Tak} the
two comonads ${\mathbb C}$ and ${\mathbb D}$ are shown to have equivalent
comodule categories.

\section{Preliminaries} \label{sec:preli}

In this section we recall some categorical preliminaries that will be used
later on.

\subsection{Notations}
Throughout the paper, in the 2-category $\mathrm{CAT}$ of
categories -- functors -- natural morphisms, the following
notations are used. The identity functor on a category ${\mathcal
C}$ is denoted by the same symbol ${\mathcal C}$. Similarly, for
any functor $F$, the identity natural morphism $F\to F$ is denoted
by the same symbol $F$. We denote vertical composition by $\circ$
and horizontal composition by juxtaposition. For
three parallel functors $F,F',F'':{\mathcal C}\to {\mathcal D}$
and natural morphisms $\alpha: F\to F'$ and $\beta: F'\to F''$,
for the composite natural morphism $F\to F''$ we write $\beta\circ
\alpha$. For consecutive functors $F:{\mathcal C} \to {\mathcal
D}$ and $G:{\mathcal D} \to {\mathcal E}$, the composite functor
is denoted by $GF:{\mathcal C} \to {\mathcal E}$. Moreover, for
functors $F,F':{\mathcal C} \to {\mathcal D}$ and $G,G':{\mathcal
D} \to {\mathcal E}$, and natural morphisms $\alpha:F \to F'$,
$\beta:G \to G'$, the Godement product $(F' \beta) \circ (\alpha
G)= (\alpha G')\circ (F\beta): FG \to F'G'$ is denoted simply by
$\alpha \beta: FG \to F'G'$. For a natural morphism
$\alpha:F \to F'$, between functors $F,F':{\mathcal C} \to
{\mathcal D}$, we denote by $\alpha X$ the morphism in ${\mathcal
D}$ obtained by evaluating $\alpha$ at an object $X$ of ${\mathcal
C}$.

The vertical category of $\mathrm{CAT}$, i.e. the category of functors --
natural morphisms, is denoted by $\mathrm{Fun}$.

In any category ${\mathcal K}$, the equalizer (resp. coequalizer) of two
parallel morphisms $f$ and $g$ (if it exists) is denoted by
$\mathrm{Equ}_{\mathcal K}(f,g)$ (resp. $\mathrm{Coequ}_{\mathcal K}(f,g)$).

\subsection{Equalizers in functor categories}
We study equalizers in the category of functors ${\mathcal C}\to {\mathcal
K}$, where ${\mathcal K}$ is assumed to have equalizers.

\begin{lemma} \label{Lem: Equ}
Let $\mathcal{C}$ and $\mathcal{K}$ be categories, let
$G$,$G':\mathcal{C}\rightarrow \mathcal{K}$ be functors and $\gamma $, $\theta
:G\rightarrow G'$ be natural morphisms.
If, for every $X\in \mathcal{C}$, there exists $\mathrm{Equ}_{\mathcal{K}}(\gamma
{X},\theta {X})$, then there exists the equalizer
$(E,i)=\mathrm{Equ}_{\mathrm{Fun}}(\gamma,\theta)$ in the category of
functors. Moreover, for any object $X$ in $\mathcal{C}$, $(EX,iX)=
\mathrm{Equ}_{\mathcal{K}}(\gamma {X},\theta {X})$.
\end{lemma}

\begin{proof}
Define a functor $E:{\mathcal C} \to {\mathcal K}$ with object map $(EX,i_X)=
\mathrm{Equ}_{\mathcal{K}}(\gamma {X},\theta {X})$. For a morphism $f:X \to X'$
in ${\mathcal C}$, naturality of $\gamma$ and $\theta$ implies that
$(Gf)\circ i_X$ equalizes the parallel morphisms $\gamma X'$ and $\theta
X'$. In light of this fact, $Ef$ is defined as the unique
morphism in ${\mathcal K}$ such
that $i_{X'} \circ (Ef) = (Gf) \circ i_X$. By construction, $i$ is a natural
transformation $E \to G$ such that $\gamma \circ i = \theta \circ i$. It
remains to prove universality of $i$. Let $H:{\mathcal C} \to {\mathcal K}$ be
a functor and $\chi:H \to G$ be
a natural morphism such that $\gamma \circ
\chi =\theta \circ \chi$. Then, for any object $X$ in ${\mathcal C}$, $(\gamma
X) \circ (\chi X)=(\theta X)\circ (\chi X)$. Since
$(EX,iX)=\mathrm{Equ}_{\mathcal K} (\gamma X, \theta X)$, there is a unique
morphism $\xi_X : H X \to EX$ such that $(iX)\circ \xi_X = \chi X$. The proof
is completed by proving naturality of $\xi_X$ in $X$. Take a morphism $f:X
\to X'$ in ${\mathcal C}$. Since $i$ and $\chi$ are natural,
\begin{eqnarray*}
(iX') \circ \xi_{X'} \circ (Hf)&=&(\chi X') \circ (Hf) = (G f) \circ (\chi
  X)\\
&=& (G f) \circ (i X) \circ \xi_X = (i X') \circ (E f) \circ \xi_X.
\end{eqnarray*}
Since $i X'$ is a monomorphism, this proves naturality of $\xi$.
\end{proof}

In the case when $\mathcal{K}$ has coequalizers, the opposite category
$\mathcal{K}^{op}$ has equalizers. Thus applying Lemma \ref{Lem: Equ} to the
opposite functors ${\mathcal C}^{op}\to \mathcal{K}^{op}$,
we conclude that coequalizers of natural transformations
between functors of codomain $\mathcal{K}$ exist, and can be computed
`objectwise'.

As an immediate consequence of Lemma \ref{Lem: Equ}, we obtain

\begin{lemma}\label{Lem: Funct1}
Let $G,G':\mathcal{C}\rightarrow \mathcal{K}$ be functors, and let $\gamma
,\theta :G\rightarrow G'$ be natural morphisms.
Assume that every pair of parallel morphisms in $\mathcal{K}$ has an
equalizer and let $\left( E,i\right) =\mathrm{Equ}_{\mathrm{Fun}}\left( \gamma
,\theta
\right) $. Under these assumptions, for any functor $P:\mathcal{D}
\rightarrow \mathcal{C}$, $\mathrm{Equ}_{\mathrm{Fun}}\left( \gamma P,\theta
P\right) = (EP,iP)$.
\end{lemma}

\begin{lemma} \label{lem:split}
For any adjunction $(N,R)$, with unit $\eta$ and counit $\epsilon$,
the following diagrams are split equalizers in the
  category of functors.
\begin{eqnarray*}
&(1)&\xymatrix{ { R}\ar[rrr]^-{\eta R} &&& {
      R N R} \ar@<2pt>[rrr]^-{\eta R N R} \ar@<-2pt>[rrr]_-{R
      N \eta R} &&& { R N R N R} };\\
&(2)&\xymatrix{ { N}\ar[rrr]^-{N \eta} &&& {
      N R N} \ar@<2pt>[rrr]^-{N\eta R N} \ar@<-2pt>[rrr]_-{N
      R N \eta} &&& { N R N R N} }.
\end{eqnarray*}
\end{lemma}

\begin{proof}
By naturality, $(\eta R N) \circ \eta = (R N  \eta) \circ
\eta$,
so both diagrams are commutative forks.

(1) A natural morphism $f:F \to R N R$, such that $(\eta R
N R) \circ f = (R N  \eta R) \circ f$, factorizes uniquely through $\eta
R$ and the morphism $(R \epsilon)\circ f$. Thus (1) is an equalizer. It
is split by the morphism $R N R \epsilon$, i.e. the identities
\begin{eqnarray*}
&&(R N R \epsilon)\circ (R N  \eta R)= R N R
\qquad \textrm{and}\\
&&(R N \eta R)\circ (R N R \epsilon)\circ  (\eta R N R)=
(\eta R N R)\circ (R N R \epsilon)\circ  (\eta R N R)
\end{eqnarray*}
hold.

(2) A natural morphism $f:F \to N R N $, such that $(N \eta R
N) \circ f = (N R N  \eta) \circ f$, factorizes uniquely through $N \eta
$ and the morphism $(\epsilon N)\circ f$. Thus (2) is an equalizer. It
is split by $\epsilon NRN$.
\end{proof}

Consider an adjoint pair of functors
$(N:{\mathcal K}\to {\mathcal C}, R:{\mathcal C}\to {\mathcal K})$ with unit
$\eta$ and counit $\epsilon$.
Since split equalizers are preserved by any functor, cf. \cite[p 110
  Proposition 2]{TTT}, for any categories ${\mathcal D}$ and
${\mathcal D}'$ and functors $P:{\mathcal C} \to {\mathcal D}$ and $Q:
{\mathcal K}\to {\mathcal D}'$,
$$
Q\eta R=\mathrm{Equ}_{\mathrm{Fun}}(Q\eta R N R, Q R N \eta R)\qquad
\textrm{and}\qquad
P N \eta = \mathrm{Equ}_{\mathrm{Fun}}(P N\eta R N, PN R N\eta).
$$

\begin{lemma}\label{lem:eta.reg}
Let ${\mathcal K}$ be a category in which all equalizers exist and let
${\mathcal C}$ be any category. Consider an adjoint pair of functors
$(N:{\mathcal K}\to {\mathcal C}, R:{\mathcal C}\to {\mathcal K})$ with unit
$\eta$ and counit $\epsilon$. Then
\begin{equation}\label{eq:A}
({\mathcal K},\eta)= \mathrm{Equ}_{\mathrm{Fun}}(RN\eta,\eta RN)
\end{equation}
if and only if $\eta$ is a regular natural monomorphism.
\end{lemma}

\begin{proof}
If \eqref{eq:A} holds then $\eta$ is obviously a regular natural
monomorphism.

Conversely, assume that $\eta$ is a regular natural monomorphism. Then we
deduce from Lemma \ref{Lem: Equ} that $\eta A$ is a regular monomorphism in
${\mathcal K}$, for any object $A\in {\mathcal K}$. By \cite[p 115 Lemma
6]{TTT} we conclude that
$$
(A,\eta A)= \mathrm{Equ}_{\mathcal K}(RN\eta A ,\eta RN A),
$$
for any object $A$ in ${\mathcal K}$. Equality \eqref{eq:A} follows by
applying Lemma \ref{Lem: Equ} again.
\end{proof}

Recall from \cite[p 111]{TTT} that a diagram like in Lemma
\ref{Lem: bohm1} is said to be {\em serially commutative} if the squares that
are bordered by parallel arrows are commutative with either {\em simultaneous}
choice of the upper or lower (or left or right) arrows.

\begin{lemma}\label{Lem: bohm1}
Consider the following serially commutative diagram in an arbitrary category
$\mathcal{K}$.
\begin{equation*}
\xymatrix{ { A} \ar[rrr]^-{i} \ar[d]_-{e} &&&
{ B} \ar@<2pt>[rrr]^-{f} \ar@<-2pt>[rrr]_-{g}
\ar[d]_-{e'}&&& { C} \ar[d]_-{e''}\\ { A'}
\ar[rrr]^-{i'} \ar@<2pt>[d]_-{n\,} \ar@<-2pt>[d]^-{\,m} &&&
{ B'} \ar@<2pt>[rrr]^-{f'} \ar@<-2pt>[rrr]_-{g'}
\ar@<2pt>[d]_-{n'\,} \ar@<-2pt>[d]^-{\,m'} &&& { C'}
\ar@<2pt>[d]_-{n''\,} \ar@<-2pt>[d]^-{\,m''} \\
{ A''} \ar[rrr]^-{i''} &&& { B''}
\ar@<2pt>[rrr]^-{f''} \ar@<-2pt>[rrr]_-{g''} &&& { C''}
}
\end{equation*}
Assume that all columns are equalizers and also the second and third rows
are equalizers. Then the first row is an equalizer too.
\end{lemma}

\begin{proof}
In order to see that the first row is a fork, note that, by commutativity of
the diagram and fork property of the second row,
$$
e''\circ f \circ i = f'\circ e'\circ i = f'\circ i'\circ e =
g' \circ i'\circ e = g' \circ e'\circ i = e'' \circ g \circ i.
$$
Since $e''$ is a monomorphism, this proves that the first row is a fork.
Take any morphism $x:X\rightarrow B$ such that $f\circ x=g\circ x$. Then
\begin{equation*}
f^{\prime }\circ e^{\prime }\circ x=e^{\prime \prime }\circ f\circ
x=e^{\prime \prime }\circ g\circ x=g^{\prime }\circ e^{\prime }\circ x.
\end{equation*}
Since the second row is an equalizer by assumption, there is a unique
morphism $y:X\rightarrow A^{\prime }$ such that
\begin{equation}
i^{\prime }\circ y=e^{\prime }\circ x.  \label{form:i'}
\end{equation}
Then $i^{\prime \prime }\circ n\circ y=i^{\prime \prime }\circ m\circ y$ and,
since $i^{\prime \prime }$ is monic, $n\circ y=m\circ y$.
Since the first column is an
equalizer, there exists a unique morphism $z:X\rightarrow A$ such that
\begin{equation}
e\circ z=y.  \label{form:e}
\end{equation}
Then $e^{\prime }\circ i\circ z=e^{\prime }\circ x$
and since $e^{\prime }$ is monic, $i\circ z=x$.
Since $e^{\prime }\circ i=i^{\prime }\circ e$ and $e^{\prime },e,i^{\prime }$
are monic, we deduce that $i$ is monic.
\end{proof}
\begin{corollary}\label{cor:pres_eq}
Let $G$,$G':\mathcal{C}\rightarrow \mathcal{K}$ be functors and $\gamma $,
$\theta :G\rightarrow G'$ be natural morphisms.
Assume that in ${\mathcal K}$ there exists the equalizer of any parallel pair
of morphisms hence there exists
$(C,i)=\mathrm{Equ}_{\mathrm{Fun}}(\gamma,\theta)$, cf. Lemma \ref{Lem:
  Equ}. Then the functor $C$ preserves equalizers provided that $G$
and $G'$ preserve equalizers.
\end{corollary}

\begin{proof}
Consider an equalizer $(E,e)=\mathrm{Equ}_{\mathcal C}(f,g)$ of morphisms
$f,g:X \to Y$ in ${\mathcal C}$.
The following diagram (in ${\mathcal K}$) is serially commutative by
naturality.
\begin{equation*}
\xymatrix{
{ {CE}} \ar[rr]^-{_{C e}} \ar[d]_-{iE} &&
{ {CX}}
\ar@<2pt>[rrr]^-{_{C f}}
\ar@<-2pt>[rrr]_-{_{C g}} \ar[d]_-{_{iX}}&&&
{ {CY}} \ar[d]^-{_{iY}}\\
{ {G E}} \ar[rr]^-{_{G e}}
\ar@<2pt>[d]_-{\gamma E\,} \ar@<-2pt>[d]^-{\,\theta E}&&
{ {G X} }
\ar@<2pt>[rrr]^-{_{G f}}
\ar@<-2pt>[rrr]_-{_{G g}}
\ar@<2pt>[d]_-{\gamma X\,}
\ar@<-2pt>[d]^-{\, \theta X} &&&
{ { G Y}}
\ar@<2pt>[d]_-{\gamma Y\, }
\ar@<-2pt>[d]^-{\, \theta Y} \\
{ {G' E} }
\ar[rr]^-{_{G' e}}&&
{ {G' X}}
\ar@<2pt>[rrr]^-{_{G' f}}
\ar@<-2pt>[rrr]_-{_{G' g}} &&&
{ {G' Y}} }
\end{equation*}
The columns are equalizers by Lemma \ref{Lem: Equ}.
The second and third rows are equalizers by assumption.
Thus the first row is an equalizer by Lemma \ref{Lem: bohm1}.
\end{proof}

\subsection{(Co)monads and their (co)modules}
We recall some basic facts about monads and their modules, mainly to fix
notation and terminology.

\begin{definition}\label{def:monad}
(1) A \emph{monad} on a category $\mathcal{K}$ is a triple
$\mathbb{T}=(T,m,u)$
where $T:\mathcal{K}\rightarrow \mathcal{K}$ is a functor
and $m:TT\rightarrow T$, $u:{\mathcal{K}}\rightarrow T $
are natural morphisms,
called the {\em product} and {\em unit}, respectively, such that
\begin{equation*}
m\circ (T m) =m\circ (m T) \qquad \text{and}\qquad
m\circ (T u) = {T} = m\circ (u T) \text{.}
\end{equation*}

(2) A \emph{morphism} between two monads
$\mathbb{T}=(T,m,u)$ and $\mathbb{T}^{\prime }=(T^{^{\prime
}},m^{\prime },u^{\prime })$ on the same category ${\mathcal K}$
is a natural morphism $\varphi :T\rightarrow T^{\prime }$ such that
\begin{equation*}
\varphi \circ m =m^{\prime }\circ (\varphi \varphi)\qquad \text{and}
\qquad \varphi \circ u=u^{{\prime }}.
\end{equation*}

(3) A $\mathbb{T}$\emph{-module} over a monad $\mathbb{T}=(T,m,u)$ on
$\mathcal{K}$ is a pair $\left( M,\mu \right)$ where $M$ is an object and
$\mu:TM \to M$ is a morphism in $\mathcal{K}$ such that
\begin{equation*}
\mu \circ (m M)=\mu \circ (T \mu ) \qquad
\text{and}\qquad \mu \circ (u M)=M \text{.}
\end{equation*}

(4) A \emph{morphism} between two $\mathbb{T}$-modules
$( M,\mu) $ and $\left( M^{\prime },\mu' \right) $
is a morphism $f:M\rightarrow M^{\prime }$ in $\mathcal{K}$ such that
$$
\mu' \circ (T f) =f\circ \mu.
$$
We denote by $\mathcal{K}_{\mathbb{T}}$ the category of
$\mathbb{T}$-modules and their morphisms.

(5) Corresponding to a monad $\mathbb{T}=(T,m,u)$ on $\mathcal{K}$, there is
an adjunction
$$
F_{\mathbb{T}}:\mathcal{K}\rightarrow \mathcal{K}_{\mathbb{T}} \qquad
U_{\mathbb{T}}:\mathcal{K}_{\mathbb{T}}\rightarrow \mathcal{K},
$$
where $U_{\mathbb{T}}$ is the {\em forgetful functor}, with object map
$(M,\mu ) \mapsto M$
and acting on the morphisms as the identity map.
$F_{\mathbb{T}}$ is the so called \emph{free functor}, with object map
$X\mapsto (TX,m X)$
and acting on the morphisms as $f\mapsto Tf$.
Note that $U_{\mathbb{T}}F_{\mathbb{T}}=T$. The unit of the adjunction is
given by
$$
u : {\mathcal K} \to U_{\mathbb{T}} F_{\mathbb{T}} =T.
$$
For a ${\mathbb{T}}$-module
$(M,\mu)$,
the counit of the adjunction
is given by
$$
 \mu:F_{\mathbb{T}}U_{\mathbb{T}} (M,\mu) =
(TM,m) \to (M, \mu).
$$
We will use the notation $\lambda^T$ for the natural transformation
$:F_{\mathbb{T}}U_{\mathbb{T}} \to {\mathcal K}_{\mathbb T}$, for which
$U_{\mathbb T} \lambda^T(M,\mu)=\mu$.

(6) A {\em comonad} on ${\mathcal K}$ is a monad on the opposite category
${\mathcal K}^{op}$. That is, a comonad ${\mathbb C}$ consists of a functor
$C:{\mathcal K}\to {\mathcal K}$, and two natural transformations $\Delta:C
\to CC$ and $\varepsilon:C \to {\mathcal K}$, called the {\em coproduct} and
{\em counit}, respectively, subject to coassociativity and counitality
constraints. Morphisms, comodules and morphisms of comodules for a comonad are
defined as morphisms, modules and morphisms of modules, respectively, for the
corresponding monad on ${\mathcal K}^{op}$. In particular, the category of
${\mathbb C}$-comodules is denoted by ${\mathcal K}^{\mathbb C}$. The
forgetful functor $U^{\mathbb C}:{\mathcal K}^{\mathbb C}\to {\mathcal K}$ has
a {\em right} adjoint $F^{\mathbb C}$ with object and morphism maps
$$
N \mapsto (CN,\Delta N)\qquad \textrm{and}\qquad f \mapsto Cf,
$$
respectively. The unit of the adjunction is given by the coaction $\varrho:M
\to CM$, for any object $(M,\varrho)\in {\mathcal K}^{\mathbb C}$ and it will
be denoted by $\gamma^{C}$. That is, $U^{\mathbb C}
\gamma^C(M,\varrho)=\varrho$. The counit is given by $\varepsilon N:CN \to N$,
for all $N\in {\mathcal K}$.
\end{definition}

\begin{proposition} \label{prop:eq_in_K_T}
Let ${\mathcal K}$ be a category with equalizers and ${\mathbb T}=(T,m,u)$ be
a monad on
${\mathcal K}$. Then any parallel pair of morphisms in ${\mathcal K}_{\mathbb
  T}$ possesses an equalizer. Moreover, the forgetful functor $U_{\mathbb T}:
{\mathcal K}_{\mathbb T} \to {\mathcal K}$ preserves and reflects equalizers.
\end{proposition}

\begin{proof}
Consider two parallel morphisms $f,g:(X,x) \to (Y,y)$ in ${\mathcal K}_{\mathbb
  T}$ and denote $(E,i):=\mathrm{Equ}_{\mathcal K}(U_{\mathbb T} f, U_{\mathbb
  T} g)$. Since $f$ and $g$ are ${\mathbb T}$-module morphisms, $x \circ (T
i)$ equalizes $f$ and $g$. So there exists a unique morphism $e:TE \to E$ in
${\mathcal K}$, such that $i\circ e = x \circ (Ti)$. By associativity and
unitality of the ${\mathbb T}$-action $x$,
\begin{eqnarray*}
i\circ e\circ (Te) &=& x\circ (Ti)\circ (Te) = x\circ (Tx) \circ (TTi) =
x\circ (m X)\circ (TT i) \\
&=& x\circ (T i) \circ (m E) = i\circ e \circ (m
E),\\
i\circ e \circ (u E) &=& x\circ (T i)\circ (u E)  = x\circ (u X) \circ i
=i.
\end{eqnarray*}
Since $i$ is monic, we conclude that $(E,e)$ is a ${\mathbb T}$-module and $i$
lifts to a ${\mathbb T}$-module morphism $\widehat i$, such that $f\circ
{\widehat i}= g\circ {\widehat i}$. It remains to prove universality of
${\widehat i}$. Consider a ${\mathbb T}$-module morphism $h:(Z,z)\to (X,x)$,
such that
$f\circ h = g\circ h$. Then $(U_{\mathbb T}f)\circ (U_{\mathbb T}h) =
(U_{\mathbb T}g) \circ (U_{\mathbb T}h)$, so by universality of $i$, there
exists a unique morphism $k:Z \to E$ in ${\mathcal K}$, such that $i \circ k =
U_{\mathbb T} h$. Since $h$ is a ${\mathbb T}$-module morphism,
$$
i\circ e \circ (Tk)= x\circ (Ti) \circ (Tk)= x\circ (T U_{\mathbb T} h) =
(U_{\mathbb T} h) \circ z = i\circ k \circ z.
$$
Since $i$ is monic, this proves that $k$ lifts to a ${\mathbb T}$-module
morphism $(Z,z) \to (E,e)$.
Since $U_{\mathbb T}$ is a right adjoint, it preserves equalizers. (Note
that it is manifest also by the above construction of
$\widehat i$ that $U_{\mathbb T}$ preserves equalizers.)
It follows by the faithfulness of $U_{\mathbb T}$ that it also reflects
equalizers. 
\end{proof}

From Proposition \ref{prop:eq_in_K_T} and Lemma \ref{Lem: Equ}, we obtain

\begin{corollary}
Let $\mathcal{K}$ be a category with equalizers and let
$\mathbb{T}=(T,m,u)$ be a monad on $\mathcal{K}$. Let $G$,$G':
\mathcal{C}\rightarrow \mathcal{K}_{\mathbb{T}}$ be functors, and let $
\gamma $, $\theta :G\rightarrow G'$ be natural morphisms. Then there
exists $\mathrm{Equ}_{\mathrm{Fun}}\left( \gamma ,\theta \right) $ and $
U_{\mathbb{T}}\mathrm{Equ}_{\mathrm{Fun}}\left( \gamma ,\theta \right)
=\mathrm{Equ}_{\mathrm{Fun}}\left( U_{\mathbb{T}}\gamma,U_{\mathbb{T}}\theta
\right) $.
\end{corollary}

\begin{proof} By Lemma \ref{Lem: Equ}, there exists $(E', \iota'):=
  \mathrm{Equ}_{\mathrm{Fun}} \left(U_{\mathbb{T}} \gamma,U_{\mathbb{T}} \theta
  \right)$. By Proposition  \ref{prop:eq_in_K_T} and Lemma \ref{Lem: Equ},
  there exists $(E,\iota):=\mathrm{Equ}_{\mathrm{Fun}} \left(\gamma,\theta
  \right)$.
  Moreover, for any object $X$ in ${\mathcal C}$,
$$
\iota' X = \mathrm{Equ}_{\mathrm{Fun}}\left(U_{\mathbb{T}} \gamma
X,U_{\mathbb{T}}
\theta X \right) = U_{\mathbb{T}}\mathrm{Equ}_{\mathrm{Fun}}\left(\gamma
X,\theta X \right) =
U_{\mathbb{T}} \iota X,
$$
where in the second equality we used that $U_{\mathbb{T}}$ preserves
equalizers. Thus the claim follows by Lemma \ref{Lem: Equ}.
\end{proof}

A statement of similar generality does not hold for coequalizers. Instead, we
have the following

\begin{proposition}\label{prop: coequal}
Let $\mathcal{K}$ be a category with coequalizers and
$\mathbb{T}=(T,m,u)$ be a monad on $\mathcal{K}$ such that $T$
preserves coequalizers. Then every parallel pair of morphisms has a
coequalizer in $\mathcal{K}_{\mathbb{T}}$. Moreover, $U_{\mathbb{T}}$
preserves and reflects coequalizers.
\end{proposition}

\begin{proof}
For two parallel morphisms $f,g:(Y,y) \to (Z,z)$ in ${\mathcal K}_{\mathbb
  T}$, denote $\left( P,\pi \right)=$ $\mathrm{Coequ}_{\mathcal{K}}\left(
U_{\mathbb{T}}
f,U_{\mathbb{T}}g\right) $. Since $f$ and $g$ are ${\mathbb T}$-module
morphisms, the morphism $\pi \circ z$ coequalizes
  $TU_{\mathbb{T}} f$ and $TU_{\mathbb{T}} g$. Therefore, by universality of
  the coequalizer $(TP,T\pi)$, there exists a unique morphism
  $\nu:TP\rightarrow P$ such that
\begin{equation}\label{eq:nu_form}
\nu \circ (T\pi) =\pi \circ z.
\end{equation}
By associativity of $z$,
$\nu \circ (T\nu ) \circ (TT\pi)= \nu \circ (m P)\circ (T
T\pi)$.
Since $T$ preserves coequalizers, $TT\pi$ is an epimorphism so we conclude
that $\nu$ is an associative $T$-action on $P$. Similarly, by unitality of
$z$, $\nu \circ (u P)\circ \pi = \pi$, hence the action $\nu$ is also
unital. Therefore, there is a ${\mathbb T}$-module morphism ${\widehat
\pi}: (Z,z) \to (P,\nu)$, such that  $U_{\mathbb{T}} {\widehat \pi}=\pi$,
cf. \eqref{eq:nu_form}.

We claim that $((P,\nu),{\widehat \pi})=\mathrm{Coequ}_{\mathcal{K}_{\mathbb
T}}(f,g)$.
Let $h:\left( Z,z\right) \rightarrow \left( W,w\right) $ be a
morphism in $\mathcal{K}_{\mathbb{T}}$ such that $h\circ f=h\circ g.$ Then
there exists a unique morphism $t:P \rightarrow W$ in $\mathcal{K}$ such
that $t\circ \pi =U_{\mathbb{T}}h$. Since $h\mathcal{\ }$is a
morphism in $\mathcal{K}_{\mathbb{T}}$,
\begin{equation*}
t\circ \nu \circ (T\pi) = w\circ (T t) \circ (T\pi).
\end{equation*}
Thus, since $T\pi$ is epi, $t$ gives rise to a morphism ${\widehat t}:\left(
P,\nu \right) \rightarrow \left( W,w\right) .$ From the uniqueness of $t$
in $\mathcal{K}$ one obviously gets the uniqueness of $\widehat{t}$ in
$\mathcal{K}_{\mathbb{T}}$. It is clear from the above construction that
$U_{\mathbb T}$ preserves coequalizers, so we conclude by uniqueness of a
coequalizer and faithfulness of $U_{\mathbb T}$ that $U_{\mathbb T}$ also
reflects coequalizers.
\end{proof}

Proposition \ref{prop: coequal} has the following immediate
consequence. Consider a category $\mathcal{K}$ with coequalizers and a monad
$\mathbb{T}$ on $\mathcal{K}$ such that the underlying functor preserves
coequalizers. Let $G$ and $G'$ be functors from any category ${\mathcal C}$ to
${\mathcal K}_{\mathbb T}$ and $\varphi,\gamma:G \to G'$ be natural
morphisms. By Proposition \ref{prop: coequal}, for any object $Z$ in
${\mathcal C}$, there exists $(K_Z,\pi_Z):=\mathrm{Coequ}_{{\mathcal K}_{\mathbb
T}}(\varphi Z,\gamma Z)$. It follows by the dual form of Lemma \ref{Lem: Equ}
that this construction defines a functor $K:{\mathcal C} \to {\mathcal
  K}_{\mathbb T}$ and a natural morphism $\pi:G' \to K$, such that
$(K,\pi)=\mathrm{Coequ}_{\mathrm{Fun}}(\varphi,\gamma)$.

\begin{definition}\label{def:comod_funct}
(1)
A {\em left comodule functor} for a comonad
${\mathbb C}=(C,\Delta,\varepsilon)$ on a category
${\mathcal A}$ is a pair $(L,l)$, where $L:{\mathcal B}\to {\mathcal A}$ is a
functor from any category ${\mathcal B}$ to ${\mathcal A}$ and $l:L\to CL$ is
a natural morphism (called a {\em coaction}) satisfying the counitality and
coassociativity constraints
$$
(\varepsilon L)\circ l = L
\qquad \textrm{and}\qquad
(\Delta L) \circ l = (Cl)\circ l.
$$

(2)
Symmetrically, a {\em right ${\mathbb C}$-comodule functor} is a pair $(R,r)$,
where $R:{\mathcal A}\to {\mathcal B}$ is a
functor from ${\mathcal A}$ to any category ${\mathcal B}$ and $r:R\to RC$ is
a natural morphism satisfying the counitality and coassociativity
constraints
$$
(R \varepsilon )\circ r = R
\qquad \textrm{and}\qquad
(R \Delta ) \circ r = (r C)\circ r.
$$

(3)
For two comonads ${\mathbb C}$ on ${\mathcal A}$ and ${\mathbb D}$ on ${\mathcal
B}$, a {\em bicomodule functor} is a triple $(Q,l,r)$, where $Q:{\mathcal
B}\to {\mathcal A}$ is a functor, $l$ is a left ${\mathbb C}$-coaction and
$r$ is a right ${\mathbb D}$-coaction on $Q$, such that
$$
(Cr)\circ l = (l D)\circ r.
$$

(4)
{\em Module functors} for a monad are defined as comodule functors for the
corresponding comonad on the opposite category.
\end{definition}

\begin{theorem}\label{theo:adjunct}
Let $\mathbb{A}=(A,m^{A},u^{A})$ and $\mathbb{T}=(T,m^{T},u^{T})$ be monads on
a category $\mathcal{K}$ and $\alpha :\mathbb{A}\rightarrow \mathbb{T}$ be a
morphism of monads. Then there exists a functor $R_A:{\mathcal K}_{\mathbb
T}\to {\mathcal K}_{\mathbb A}$ such that
\begin{equation}
U_{\mathbb{A}}R_{A}=U_{\mathbb{T}}\qquad \text{and}\qquad
U_{\mathbb{A}} \lambda^A R_A=(U_{\mathbb{T}}\lambda^T) \circ (\alpha
{U_{\mathbb{T}}}).  \label{form:UTdiphistar}
\end{equation}
Moreover, if there exist coequalizers in ${\mathcal K}$ and $T$ preserves
coequalizers, then $R_A$ has a left adjoint.
\end{theorem}

\begin{proof}
In order to construct a functor $R_A$, note that
for $(X,x)\in \mathcal{K}_{\mathbb{T}}$, $(X,x\circ (\alpha X)) \in
\mathcal{K}_{\mathbb{A}}$. A morphism $ f:(X,x)\rightarrow
(X^{\prime},x^{\prime })$ in $\mathcal{K}_{\mathbb{T}}$ can be regarded as a
morphism $(X,(x\circ \alpha X))\rightarrow (X^{\prime },x^{\prime }\circ
(\alpha X^{\prime }))$ in $\mathcal{K}_{\mathbb{A}}$. Therefore we
introduce the functor $R_{A}:\mathcal{K}_{\mathbb{T}}$ $\rightarrow \mathcal{
K}_{\mathbb{A}}$ by setting \eqref{form:UTdiphistar}.

Assume now that the category ${\mathcal K}$ has coequalizers.
We can define a functor $N_A:{\mathcal K}_{\mathbb A} \to {\mathcal
K}_{\mathbb T}$ as an `${\mathbb A}$-module product' of the right ${\mathbb
  A}$-module functor $(F_{\mathbb T}, (\lambda^T F_{\mathbb T})\circ
(F_{\mathbb T} \alpha))$ and the left ${\mathbb A}$-module functor
$(U_{\mathbb A}, U_{\mathbb A} \lambda^A)$.
That is, we define $N_A$ via the coequalizer
\begin{equation}\label{eq:N_A}
\left( N_{A},\chi ^{A}\right) =\mathrm{Coequ}_{\mathrm{Fun}}\left(
(\lambda^T {F_{\mathbb{T}}U_{\mathbb{A}}})\circ (F_{\mathbb{T}} \alpha
{U_{\mathbb{A}}}),F_{\mathbb{T}} U_{\mathbb{A}} \lambda ^{{A}}\right).
\end{equation}
It takes a morphism
$f:\left( Y,y\right) \rightarrow \left( Y^{\prime
},y^{\prime }\right) $ in
$\mathcal{K}_{\mathbb{A}}$ to the the unique morphism $N_{A} f$
in $\mathcal{K}_{\mathbb{T}}$ for which
$$
(N_{A}f) \circ (\chi^A ( Y,y) )=(\chi ^A \left(
Y^{\prime },y^{\prime }\right))\circ (F_{\mathbb{T}}U_{\mathbb{A}}f) .
$$
Under the assumption that $T$ (hence by Proposition \ref{prop: coequal} also
$U_{\mathbb T}$) preserves coequalizers, we prove next
that $N_{A}$ and $R_{A}$ are adjoint functors. By the
construction of $R_A$, the morphism $\lambda^T$ coequalizes
the parallel morphisms $(\lambda^T F_{\mathbb T} {U_{\mathbb{A}}R_{A}})\circ
(F_{\mathbb{T}}\alpha {U_{\mathbb{A}}R_{A}})$ and
$F_{\mathbb{T}}U_{\mathbb{A}}\lambda^A R_{A}$. Thus it follows by universality of
the coequalizer $\chi^A{R_{A}}$ (cf. dual form of Lemma \ref{Lem: Funct1})
that there exists a unique natural morphism $\epsilon
^{A}:N_{A}R_{A}\rightarrow {\mathcal{K}_{\mathbb{T}}}$ such that
\begin{equation}
\epsilon ^{A}\circ (\chi^A {R_{A}})=\lambda ^T\text{.}
\label{form: epsil}
\end{equation}
Since $(U_{\mathbb{T}} \chi ^{A}) \circ (u^T {U_{\mathbb{A}}})$ is an
${\mathbb A}$-module morphism in the sense that
$$
(U_{\mathbb T} \lambda^T N_A) \circ (\alpha {U_{\mathbb{T}}N_{A}}) \circ (A
U_{\mathbb{T}} \chi ^{A}) \circ (A u^T {U_{\mathbb{A}}}) =
(U_{\mathbb{T}} \chi ^{A}) \circ (u^T {U_{\mathbb{A}}}) \circ (U_{\mathbb A}
\lambda ^A),
$$
it lifts to a natural morphism
$\eta^{A}:{\mathcal{K}_{\mathbb{A}}}\rightarrow R_{A}N_{A}$ in the sense that
\begin{equation}
U_{\mathbb{A}} \eta ^{A} = (U_{\mathbb{T}} \chi ^{A})
\circ (u^T {U_{\mathbb{A}}}) . \label{form:UAdietaA}
\end{equation}
Note that \eqref{form:UAdietaA} immediately implies
\begin{equation}
(\lambda^T {N_{A}}) \circ (F_{\mathbb{T}}U_{\mathbb{A}}\eta^{A}) =\chi ^{A}.
  \label{form:agg0,5}
\end{equation}
From \eqref{form: epsil} and \eqref{form:agg0,5} we deduce that
$
(U_{\mathbb T}\epsilon^A {N_{A}}) \circ (U_{\mathbb T} N_{A} \eta ^{A})\circ
(U_{\mathbb T} \chi ^{A})=
U_{\mathbb T}\chi ^{A}.
$
Since $U_{\mathbb T}\chi^A$ is epi and $U_{\mathbb T}$ is faithful, this implies
$(\epsilon^A {N_{A}})\circ (N_{A}\eta ^{A})= {N_{A}}$.
Similarly, $(U_{\mathbb{A}} R_{A} \epsilon ^{A})\circ (U_{\mathbb{A}}\eta^A
{R_{A}})= U_{\mathbb{A}}{R_{A}}$, so by faithfulness of $U_{\mathbb A}$,
$(R_{A}\epsilon ^{A})\circ (\eta^A {R_{A}})={R_{A}}$. Thus we proved that
$(N_{A},R_{A})$ is an adjunction.
\end{proof}

Theorem \ref{theo:adjunct} implies, in particular, that
$(N_{A}R_{A},N_{A}\eta^A {R_{A}},\epsilon ^{A})$ is a
comonad on $\mathcal{K}_{\mathbb{T}}$ and $(R_{A}N_{A},R_{A}\epsilon^A
{N_{A}},\eta ^{A})$ is a monad on $\mathcal{K}_{\mathbb{A}}$.

For the adjoint functors in Theorem \ref{theo:adjunct}, induced by a monad
morphism $\alpha$, we use also the notation $N_A=\alpha^*$ and $R_A=\alpha_*$.
Note that, for two monad morphisms $\alpha :\mathbb{A}\rightarrow \mathbb{T}$
and $\varphi: \mathbb{T} \to \mathbb{T}'$, $(\varphi\circ \alpha)_\ast =
\alpha_\ast\circ \varphi_\ast$ and $(\varphi\circ \alpha)^\ast =\varphi^\ast
\circ  \alpha^\ast$.
The identity functor on any category $\mathcal{K}$ is a monad, via the
multiplication and unit given by the identity natural morphism
$\mathcal{K}$.
Given any monad $\mathbb{T=}\left( T,m,u\right) $ on $\mathcal{K}$ ,
$u$ gives rise to a monad morphism  ${\mathcal{K}}\rightarrow \mathbb{T}$
and $U_{\mathbb{T}}=u_{\ast }$
while $F_{\mathbb{T}}=u^{\ast }$.

\begin{theorem}\label{theo:main}
Let $\mathcal{K}$ be a category with coequalizers
and let $\alpha :\mathbb{A}=(A,m^{A},u^{A}) \to \mathbb{T}=(T,m^{T},u^{T})$ be
a morphism of monads on $\mathcal{K}$.
Assume that $A$ and $T$ preserve coequalizers.
Consider the canonical adjunction $(N_A,R_A)$, associated to $\alpha$ in
Theorem \ref{theo:adjunct}. Then the Eilenberg-Moore comparison functor
$K:\mathcal{K}_{\mathbb{T}}\rightarrow \left(\mathcal{K}_{\mathbb{A}}\right)
_{{R_A N_A}}$ is an isomorphism.
\end{theorem}

\begin{proof}
Denote by $U$ forgetful functor $\left(\mathcal{K}_{\mathbb{A}}\right) _{R_A
  N_A}\rightarrow \mathcal{K}_{\mathbb{A}}$. Recall that the comparison
  functor $K$ has the explicit form
\begin{equation*}
K\left( X,x\right) =\left( R_{A}\left( X,x\right) ,R_{A}
\epsilon^A {\left( X,x\right) } \right) =\left( (X,x\circ (\alpha X)),
R_{A}\epsilon^A {\left( X,x\right)} \right),
\end{equation*}
for any object $(X,x)$ of ${\mathcal K}_{\mathbb T}$,
and $K f=R_{A}f$, for each morphism $f$ in $\mathcal{K}_{\mathbb{T}}$.
In what follows we construct the inverse ${\widetilde K}$ of $K$. A direct
computation shows that, for any object $((X,x), x')\in ({\mathcal
  K}_A)_{R_AN_A}$, the morphism
$(U_{\mathbb{A}} x') \circ (U_{\mathbb{T}}\chi^A ( X,x))$ is a ${\mathbb
  T}$-action on $X$. Moreover, with respect to this ${\mathbb T}$-action, any
morphism in $({\mathcal K}_A)_{R_A N_A}$ is a morphism of ${\mathbb
  T}$-modules.
Thus, for any object $((X,x),x')\in ({\mathcal K}_A)_{R_AN_A}$ we can put
\begin{equation*}
\widetilde{K}((X,x),x')=(X,(U_{\mathbb{A}}x') \circ
(U_{\mathbb{T}}\chi^A {\left( X,x\right) }) ).
\end{equation*}
For a morphism $f$ in $({\mathcal K}_A)_{R_AN_A}$, $\widetilde{K}f$ is
defined as the unique morphism in $ \mathcal{K}_{\mathbb{T}}$ such that
$U_{\mathbb{T}}\widetilde{K} f=U_{\mathbb{A}}Uf$. The equality
${\widetilde K} K= {\mathcal K}_{\mathbb T}$ is obvious. The proof of $K
{\widetilde K}=({\mathcal K}_A)_{R_A N_A}$ is slightly longer but
straightforward, so it is left to the reader too.
\end{proof}

\begin{definition}\label{def:entw}
Let $\mathbb{A}=(A,m,u)$ be a monad and $\mathbb{C}=(C,\Delta,\varepsilon)$ be
a comonad on the same category $\mathcal{K}$. A
natural morphism $\Psi :AC\rightarrow CA$ is called
a {\em mixed distributive law} (or in some papers an \emph{entwining}) if
\begin{itemize}
\item $\Psi \circ (m {C}) =( C m) \circ (\Psi {A})\circ(A\Psi)$ and
$\Psi \circ (u C) =Cu$,
\item $(\Delta A) \circ \Psi =( C \Psi) \circ (\Psi {C}) \circ (A\Delta) \, $
and $(\varepsilon A) \circ \Psi =A\varepsilon $.
\end{itemize}
\end{definition}

\begin{theorem}\label{lift}\cite{Be}
Let $\mathbb{A}=(A,m,u)$ be a monad and $\mathbb{C}=(C,\Delta,\varepsilon)$
be a comonad on a category $\mathcal{K}$.
There is a bijection between
\begin{itemize}
\item {\em liftings} of ${\mathbb C}$ to a comonad $\widetilde{\mathbb C}$ on
  ${\mathcal K}_{\mathbb A}$, i.e. comonads
$\widetilde{\mathbb{C}}=({\widetilde C},{\widetilde \Delta},{\widetilde
  \varepsilon})$ on $\mathcal{K}_{\mathbb{A}}$, such that
$U_{\mathbb A} {\widetilde C} = C U_{\mathbb A}$,
$U_{\mathbb A} {\widetilde \Delta}=\Delta U_{\mathbb A}$ and  $U_{\mathbb A}
{\widetilde \varepsilon}=\varepsilon  U_{\mathbb A}$;
\item mixed distributive laws $\Psi:AC \to CA$.
\end{itemize}
\end{theorem}

\begin{proof}
It is a standard result due to Johnstone \cite{Johnst} that a mixed
distributive law $AC \to CA$ determines a functor
${\widetilde C}:\mathcal{K}_{\mathbb{A}}\rightarrow \mathcal{K}_{\mathbb{A}}$,
with object map
\begin{equation*}
{\widetilde C}(X,x)=(CX,(C x)\circ (\Psi {X}))
\end{equation*}
and morphism map satisfying $U_{\mathbb{A}} {\widetilde C}f=C
U_{\mathbb{A}} f$, for every morphism $f\in
\mathcal{K}_{\mathbb{A}}$. It is proven by a slightly twisted
version of Beck's arguments \cite{Be} (cf. \cite[5.1]{W}) that it is a
comonad, with the stated coproduct and counit.

Conversely, if $\widetilde{\mathbb C}$ is a lifting of ${\mathbb C}$ then a
mixed distributive law is constructed as
$$
\xymatrix{
AC \ar[r]^-{AC u}&A C A=U_{\mathbb A} F_{\mathbb A} U_{\mathbb A} {\widetilde C}
F_{\mathbb A} \ar[rr]^-{U_{\mathbb A} \lambda^A {\widetilde C} F_{\mathbb A}}&&
U_{\mathbb A} {\widetilde C} F_{\mathbb A} = C A.
}
$$
\end{proof}

\section{Galois functors for mixed distributive laws}
\label{sec:entw}

In this section we study particular kinds of comonads -- those arising from
mixed distributive laws -- including the comonads arising from entwining
structures of algebras and coalgebras. Our aim is to reformulate the Galois
property of a comodule functor of such a comonad, in terms of so called
regular comonad arrows. Comonad arrows (not only those corresponding to Galois
functors for mixed distributive laws) will play an important role in later
sections: They will be related to pre-torsors.
Together with the results of the current section, this implies a
relation between Galois functors and pre-torsors.

\subsection{Galois functors and regular comonad arrows}\label{sec:com.tr}

Under various names, functors that we term {\em Galois functors} (following
\cite[Definition 4.5]{MW}), have been discussed by several authors, see e.g.
\cite{Dubuc}, \cite{KeSt}, \cite{GT}, \cite{MW}.

\begin{definition}
If, for a left ${\mathbb C}$-comodule functor $(L,l)$, the underlying functor
$L$ has a right adjoint $R$, then there is a canonical comonad morphism
\begin{equation}\label{eq:L.can}
\mathrm{can}:= (C\epsilon)\circ (lR):LR\to C,
\end{equation}
where $\epsilon$ denotes the counit of the adjunction,
see \cite[Proposition 3.3]{KeSt}.
A {\em ${\mathbb C}$-Galois functor} is, by definition, a left
${\mathbb C}$-comodule functor $(L,l)$, such that the underlying functor $L$
possesses a right adjoint $R$ and the canonical comonad morphism
\eqref{eq:L.can} is an isomorphism.
\end{definition}

The bicategory of (co)monads was introduced in the paper \cite{S}. Its 1-cells,
recalled under the name {\em comonad arrow} in Definition
\ref{def:com.tr}, generalize morphisms of comonads.

\begin{definition}\label{def:com.tr}
Consider two comonads ${\mathbb C}=(C,\Delta,\varepsilon)$ and ${\mathbb
C}'=(C',\Delta',\varepsilon')$, on respective categories ${\mathcal A}$ and
${\mathcal A}'$.
A {\em comonad arrow} from ${\mathbb C}$ to ${\mathbb C}'$ is a pair
$(F,\xi)$, where $F:{\mathcal A}'\to {\mathcal A}$ is a functor and $\xi:CF
\to FC'$ is a natural morphism subject to the conditions
\begin{equation}\label{eq:com.tr}
(F\varepsilon')\circ \xi = \varepsilon F
\qquad \textrm{and}\qquad
(F\Delta')\circ \xi = (\xi C')\circ (C\xi)\circ (\Delta F).
\end{equation}

A comonad arrow $(F,\xi)$ is said to be {\em regular} provided that
$\xi$ is an isomorphism.

A comonad arrow $(F,\xi)$ is termed {\em co-regular} if $F$ has a left adjoint
$G$ (with unit $\eta$ and counit $\epsilon$ of the adjunction) and ${\overline
  \xi}:= (\epsilon C' G)\circ (G\xi G)\circ (GC \eta):
GC \to C' G$ is an isomorphism.
\end{definition}

In the following theorem, the functor
$$
{\mathcal A}_{{\mathbb N}'}\to {\mathcal A}_{\mathbb N},
\qquad
(X,x) \mapsto (X,x\circ (\varphi X)),
$$
induced by a morphism $\varphi:{\mathbb N}\to {\mathbb N}'$ of
monads on the category ${\mathcal A}$ (cf. Theorem \ref{theo:adjunct}), is
denoted by $\varphi_*$. For an adjunction $(L,R)$, the unit and counit are
denoted by $\eta$ and $\epsilon$, respectively. For a second adjunction
$(L',R')$, primed symbols $\eta'$ and $\epsilon'$ are used.
The symbol $\mathrm{can}$ (resp. $\mathrm{can}'$) denotes the comonad morphism
\eqref{eq:L.can} corresponding to a comodule functor $L$ (resp. $L'$).

\begin{theorem} \label{thm:Gal.com.tr}
(1)
Let ${\mathbb N}=(N,m^N,u^N)
$ be a monad and ${\mathbb C}=(C,\Delta^C,\varepsilon^C)
$ be a comonad on a category
${\mathcal A}$.
For any category ${\mathcal B}$  and any adjoint pair of functors
$(L:{\mathcal B} \to {\mathcal A}_{\mathbb N}, R:
{\mathcal A}_{\mathbb N} \to {\mathcal B})$,
there is a bijective correspondence between the sets of the following data.

(a) Liftings of ${\mathbb C}$ to a comonad ${\widetilde {\mathbb C}}$ on
${\mathcal A}_{\mathbb N}$ together with a
${\widetilde {\mathbb C}}$-Galois functor structure on $L$;

(b) regular
comonad arrows $(U_{\mathbb N},\xi)$ from ${\mathbb C}$ to the comonad $LR$.

(2) Let ${\mathbb N}$ and ${\mathbb N}'$ be monads on a
category ${\mathcal A}$. Let
${\mathbb C}$ and ${\mathbb C}'$ be comonads on ${\mathcal A}$,
such that ${\mathbb C}$ lifts to a comonad ${\widetilde {\mathbb C}}$ on
${\mathcal A}_{\mathbb N}$ and
${\mathbb C}'$ lifts to a comonad ${\widetilde {\mathbb C}'}$ on
${\mathcal A}_{{\mathbb N}'}$.
Let $(L,l)$ be a
${\widetilde {\mathbb C}}$-Galois functor with right adjoint $R$ and $(L',l')$
be a ${\widetilde {\mathbb C}'}$-Galois functor with right adjoint $R'$.
In this setting, for any monad morphism $\varphi:{\mathbb N} \to {\mathbb N}'$
and comonad morphism
$\theta:{\mathbb C} \to {\mathbb C}'$, the following groups of statements
are equivalent.
\begin{itemize}
\item[{\textit{(a)}}$\bullet$] $R\varphi_* =R'$,
\item there is a comonad arrow $(\varphi_*,{\widetilde
  \theta})$ from
${\widetilde {\mathbb C}}$ to ${\widetilde {\mathbb C}'}$, such that
  $U_{\mathbb N}{\widetilde \theta} = \theta U_{{\mathbb N}'}$,
\item the natural morphism ${\widetilde \varphi}:= (\epsilon \varphi_*
  L')\circ (L \eta'):L\to \varphi_* L'$
  satisfies
$
(\varphi_* l')\circ {\widetilde \varphi} = ({\widetilde \theta} L')\circ
  ({\widetilde {\mathbb C}} {\widetilde \varphi})\circ l
$.
\end{itemize}
\begin{itemize}
\item[{\textit{(b)}}$\bullet$] $R\varphi_* =R'$,
\item
$
(U_{{\mathbb N}'} \mathrm{can}')\circ
(U_{\mathbb N} \epsilon \varphi_* L'R')\circ(U_{\mathbb N} L \eta' R')=
  (\theta U_{{\mathbb N}'})\circ (U_{\mathbb N} \mathrm{can} \varphi_*)
$.
\end{itemize}
\end{theorem}

\begin{proof}
(1) Consider first data as in part (a) and put $\xi:=U_{\mathbb N}
{\mathrm{can}}^{-1}$. It is
obviously an isomorphism and the identities \eqref{eq:com.tr} follow by using
that ${\mathrm{can}}:LR \to {\widetilde {\mathbb C}}$ is a comonad morphism.

Conversely, in terms of the data in part (b), a mixed distributive law $NC \to
CN$ is given by the natural morphism
$$
\psi:=(\xi^{-1} F_{\mathbb N}) \circ (U_{\mathbb N} \lambda^N L R F_{\mathbb
  N}) \circ (N \xi F_{\mathbb N}) \circ (N C u^N).
$$
This proves that ${\mathbb C}$ lifts to a comonad ${\widetilde {\mathbb C}}$
on ${\mathcal A}_{\mathbb N}$, cf. Theorem \ref{lift}. Moreover, $\psi$
induces an ${\mathbb N}$-action
$$
(C U_{\mathbb N} \lambda^N) \circ (\psi U_{\mathbb N}) : NCU_{\mathbb N} \to
CU_{\mathbb N}
$$
on $CU_{\mathbb N}$. It is easy to check that (by naturality and the
adjunction relations) $\xi$ is an ${\mathbb N}$-module morphism in the sense
that 
$$
\xi\circ (C U_{\mathbb N}\lambda^N) \circ (\psi U_{\mathbb N}) =
(U_{\mathbb N} \lambda^N  L R) \circ (N \xi).
$$
This means that $\xi$ gives rise to a morphism $\widehat{\xi}$ in ${\mathcal
  A}_{\mathbb N}$, such that $U_{\mathbb N}
\widehat{\xi} =\xi$. Since the forgetful functor $U_{\mathbb N}$ reflects
isomorphisms, $\widehat{\xi}$ is an isomorphism in ${\mathcal A}_{\mathbb N}$.
This enables us to equip $L$ with a ${\widetilde {\mathbb
    C}}$-coaction by putting
\begin{equation}\label{eq:C.coac}
l:= ({\widehat \xi}^{-1} L)\circ (L \eta): L \to {\widetilde {\mathbb C}}L.
\end{equation}
Using identities \eqref{eq:com.tr}, naturality and the adjunction relations,
we find that
$$
(\varepsilon^C U_{\mathbb N} L) \circ (U_{\mathbb N} l)=U_{\mathbb N} L
\qquad\textrm{and}\qquad
(\Delta^C U_{\mathbb N} L) \circ (U_{\mathbb N} l) =
(C U_{\mathbb N} l )\circ (U_{\mathbb N} l).
$$
Since ${\widetilde {\mathbb C}}$ is the lifting of ${\mathbb C}$ and
$U_{\mathbb N}$ is faithful, this implies coassociativity and counitality of the
coaction $l$.
Moreover, by \eqref{eq:C.coac}, naturality and the adjunction relations,
$$
U_{\mathbb N} \mathrm{can}=
(C U_{\mathbb N} \epsilon)\circ (U_{\mathbb N} l R) = \xi^{-1}=
U_{\mathbb N} {\widehat \xi}^{-1}.
$$
Thus by faithfulness of $U_{\mathbb N}$, we
conclude that ${\mathrm{can}} = {\widehat \xi}^{-1}$ is an isomorphism, hence
$(L,l)$ is a ${\widetilde {\mathbb C}}$-Galois functor.

Above constructions are easily checked to yield a bijective correspondence
between the data in parts (a) and (b).

(2) Assume that the conditions in part (a) hold. By
\eqref{eq:L.can}, the second property in part (b) is equivalent to
$$
(C' U_{{\mathbb N}'} \epsilon') \circ (U_{{\mathbb N}'} l' R')\circ
(U_{{\mathbb N}} {\widetilde \varphi} R')=
(\theta U_{{\mathbb N}'}) \circ (C U_{\mathbb N} \epsilon\varphi_*) \circ
(U_{\mathbb N} l R \varphi_*),
$$
which is proven by using the third condition in part (a), the identity
$U_{\mathbb N} {\widetilde \theta} = \theta U_{{\mathbb N}'}$ and the fact
that ${\widetilde {\mathbb C}}$ is the lifting of ${\mathbb C}$, naturality
and definition of ${\widetilde \varphi}$, adjunction relations and
finally the first identity in part (a).

Conversely, assume that assertion (b) holds. The natural morphism
$\theta$ is checked to have a lifting ${\widetilde \theta}$ if and only if the
mixed distributive laws $\psi$ and $\psi'$, that are responsible for the lifting
of ${\mathbb C}$ to ${\widetilde {\mathbb C}}$ and the lifting of ${\mathbb
C}'$ to ${\widetilde {\mathbb C}}'$, respectively, satisfy
\begin{equation}\label{eq:psi.comp}
\psi' \circ (\varphi \theta) = (\theta \varphi) \circ \psi.
\end{equation}
In order to prove \eqref{eq:psi.comp}, note that both $N$ and $N'$ are left
${\mathbb 
  N}$-module functors, via the actions provided by the multiplication $m^N$ in
${\mathbb N}$, and $m^{N'}\circ (\varphi N')$, respectively. With respect to
these actions, $\varphi$ is an ${\mathbb N}$-module morphism. Hence there
exists a (unique) morphism ${\widehat \varphi}:F_{\mathbb N} \to \varphi_*
F_{\mathbb N'}$ (explicitly given by ${\widehat \varphi}=(\lambda^N \varphi_*
F_{{\mathbb N}'})\circ (F_{\mathbb N} u^{N'})$) such that $U_{\mathbb N} {\widehat
\varphi} =\varphi$.
Moreover, by the definition of the functor $\varphi_*$, the identity
$(U_{{\mathbb N}'} \lambda^{N'})\circ (\varphi U_{{\mathbb N}'})=U_{\mathbb N}
\lambda^N \varphi_*$ holds. Making use of these observations,
\eqref{eq:psi.comp} follows by
recalling the form of $\psi$ and $\psi'$ from part (1), repeated use of the
second condition in part (b), the monad morphism property of $\varphi$
and naturality. This proves the existence of ${\widetilde \theta}$.
Since $\theta$ is a comonad morphism, ${\widetilde \theta}$ satisfies
\eqref{eq:com.tr}.

The third condition in part (a)
follows by \eqref{eq:C.coac} and the analogous formula for
$l'$, the second condition in part (b), naturality and adjunction relations.
\end{proof}

Theorem \ref{thm:Gal.com.tr} can be rephrased as a statement about the
existence of a certain functor, between categories defined below.

\begin{definition} \label{def:cat_adj}
For any two categories ${\mathcal A}$ and ${\mathcal B}$,
  the category $\mathrm{Adj}({\mathcal A},{\mathcal B})$ is defined as
  follows.

{\em Objects} are triples $({\mathcal T},(N_A,R_A), (N_B,R_B))$,
  where ${\mathcal T}$ is a category and
$$
(N_A:{\mathcal A} \to {\mathcal T},R_A:{\mathcal T}\to {\mathcal A})
\quad \textrm{and}\quad
(N_B:{\mathcal B} \to {\mathcal T},R_B:{\mathcal T}\to {\mathcal B})
$$
are adjunctions. We denote the respective units of the adjunctions by $\eta^A$
and $\eta^B$ and the counits by $\epsilon^A$ and $\epsilon^B$.

{\em Morphisms} $({\mathcal T},(N_A,R_A), (N_B,R_B))\to ({\mathcal
  T}',(N'_A,R'_A), (N'_B,R'_B))$
are functors $F:{\mathcal T}' \to {\mathcal T}$ such that $R_A F = R'_A$ and
$R_B F= R'_B$.
\end{definition}

Note that a morphism $F$ in $\mathrm{Adj}({\mathcal A}, {\mathcal B})$
comes equipped with natural morphisms
$a:=(\epsilon^A F N'_A)\circ (N_A \eta^{\prime A}):N_A \to F N'_A$ and
$b:=(\epsilon^B F N'_B)\circ (N_B \eta^{\prime B}): N_B \to F N'_B$ such that
the following compatibility conditions hold.
\begin{eqnarray}\label{eq:a&b}
&(R_A a) \circ \eta^A = \eta^{\prime A}\qquad
\quad \textrm{and}\quad
&(R_B b) \circ \eta^B = \eta^{\prime B},\\
&(F \epsilon^{\prime A}) \circ (a R'_A) = \epsilon^A F
\quad \textrm{and}\quad
&(F\epsilon^{\prime B}) \circ (b R'_B) = \epsilon^B F.
\nonumber
\end{eqnarray}
In fact, $a$ and $b$ are unique natural morphisms satisfying these
identities.

An object $({\mathcal T},(N_A,R_A), (N_B,R_B))$ of $\mathrm{Adj}({\mathcal
  A},{\mathcal B})$
determines two comonads on ${\mathcal T}$,
$$
(N_A R_A, N_A \eta^A R_A,\epsilon^A)\quad \textrm{and}\quad
(N_B R_B, N_B \eta^B R_B,\epsilon^B).
$$

\begin{definition}\label{def:Gal_adj}
The category $\mathrm{Arr}({\mathcal A},{\mathcal B})$ is defined to have
{\em objects} of the form
$({\mathcal T},(N_A,R_A), (N_B,R_B), {\mathbb C},\xi)$,
where $({\mathcal T},(N_A,R_A), (N_B,R_B))$ is an object in
$\mathrm{Adj}({\mathcal A},{\mathcal B})$, ${\mathbb C}$ is a comonad on
${\mathcal A}$ and $(R_A,\xi)$ is a comonad arrow from ${\mathbb C}$
to $N_B R_B$.

A {\em morphism} in $\mathrm{Arr}({\mathcal A},{\mathcal B})$
$$
({\mathcal T},(N_A,R_A), (N_B,R_B), {\mathbb C},\xi)\to
({\mathcal T}',(N'_A,R'_A), (N'_B,R'_B), {\mathbb C}',\xi')
$$
consists of a morphism $
F:({\mathcal T},(N_A,R_A), (N_B,R_B)) \to
({\mathcal T}',(N'_A,R'_A), (N'_B,R'_B))$ in $\mathrm{Adj}({\mathcal
  A},{\mathcal B})$ and a comonad morphism $t: {\mathbb C}\to {\mathbb C}'$,
such that
\begin{equation}\label{eq:RGal_mor}
(R_A b R'_B) \circ (\xi F) = \xi'\circ (t R'_A),
\end{equation}
where the morphism $b$ was introduced in \eqref{eq:a&b}.

The full subcategory of $\mathrm{Arr}({\mathcal A},{\mathcal B})$ of
objects $({\mathcal T},(N_A,R_A), (N_B,R_B), {\mathbb C},\xi)$,
such that the comonad arrow $(R_A,\xi)$ is {\em regular}, will be denoted by
$\mathrm{RArr}({\mathcal A},{\mathcal B})$.

The full subcategory of $\mathrm{Arr}({\mathcal A},{\mathcal B})$ of
objects $({\mathcal T},(N_A,R_A), (N_B,R_B), {\mathbb C},\xi)$, such that
the comonad arrow $(R_A,\xi)$ is {\em co-regular}, will be denoted by
$\overline{\mathrm{RArr}} ({\mathcal A},{\mathcal B})$.
\end{definition}
Obviously, composites of morphisms in
$\mathrm{Arr}({\mathcal  A},{\mathcal B})$ are again morphisms
in $\mathrm{Arr}({\mathcal  A},{\mathcal B})$.

The other category occurring in Theorem \ref{thm:Gal.com.tr} is the following.

\begin{definition}
For any two categories ${\mathcal A}$ and ${\mathcal B}$,
{\em objects} of the category $\mathrm{Gal}({\mathcal A},{\mathcal B})$
are quintuples $({\mathbb N},{{\mathbb C}},\psi,L,l)$, where
${\mathbb N}$ is a monad and ${\mathbb C}$ is a comonad on
${\mathcal A}$, and $\psi$ is a mixed distributive law between them.
$L:{\mathcal B}\to {\mathcal A}_{\mathbb N}$ is a Galois functor, with
  coaction $l$,
for the lifted comonad ${\widetilde {\mathbb C}}$ on ${\mathcal
A}_{\mathbb N}$, determined by the mixed distributive law $\psi$.

A {\em morphism} $({\mathbb N},{{\mathbb C}},\psi,L,l) \to ({\mathbb
  N}',{{\mathbb C}'},\psi',L',l')$ is a pair $(\varphi,\theta)$, consisting
of a monad morphism $\varphi:{\mathbb N} \to {\mathbb N}'$ and a comonad
morphism $\theta:{\mathbb C} \to {\mathbb C}'$, subject to the conditions in
Theorem \ref{thm:Gal.com.tr} (2)(b).
\end{definition}

\begin{corollary}\label{cor:Gal.RArr}
By Theorem \ref{thm:Gal.com.tr}, there is a
functor $I: \mathrm{Gal}({\mathcal A},{\mathcal B}) \to
\mathrm{RArr}({\mathcal A},{\mathcal B})$, that is faithful and
injective also on the objects. The object and morphism maps of $I$ are
$$
({\mathbb N},{\mathbb C},\psi,L,l) \mapsto
\big({\mathcal A}_{\mathbb N}, (F_{\mathbb N},U_{\mathbb N}),(L,R),
{\mathbb C},U_{\mathbb N}
{\mathrm{can}}^{-1} \big)
\qquad \textrm{and}\qquad
(\varphi,\theta) \mapsto
(\varphi_*,\theta),
$$
where $R$ is the right adjoint of $L$.
Moreover, any object in $\mathrm{RArr}({\mathcal A},{\mathcal B})$ that is of
the form $\big({\mathcal A}_{\mathbb N}, (F_{\mathbb N},U_{\mathbb N}),(L,R),
{\mathbb C}, \xi \big)$ arises as the image of an (unique) object in
$\mathrm{Gal}({\mathcal A},{\mathcal B})$ under the functor $I$.
\end{corollary}

\subsection{Examples from bimodules} \label{sec:bimod.ex}

For an associative and unital algebra ${\bf A}$, consider an ${\bf A}$-ring
${\bf T}$ and an ${\bf A}$-coring ${\bf C}$, that are entwined by
$\boldsymbol{\psi}:{\bf C \ot_A T} \to {\bf T \ot_A C}$.
(For a review of these structures we refer to Sections A.1, A.2 and A.4 of
\cite{ABM}.)
Denote the induced
${\bf T}$-coring ${\bf T\ot_A C}$ (cf. \cite[Section A.4]{ABM}) by
${\widetilde{\bf C}}$.
These data determine a monad ${\mathbb N}= (-)\ot_{\bf A} {\bf T}$ and a
comonad ${\mathbb C}=(-)\ot_{\bf A} {\bf C}$ on the category ${\mathcal
A}:=\mathrm{Mod}$-${\bf A}$ and also a lifted comonad ${\widetilde {\mathbb
C}}=(-)\ot_{\bf T} {\widetilde{\bf C}}\cong (-)\ot_{\bf A} {\bf C}$ on
${\mathcal A}_{\mathbb N} \cong \mathrm{Mod}$-${\bf T}$.

Take now a right ${\widetilde {\bf C}}$-comodule (i.e. entwined module)
$\Sigma$, and let ${\bf B}$ be any subalgebra of $\mathrm{End}^{\widetilde{\bf
 C}}(\Sigma)$. Then
the functor $(-)\ot_{\bf B} \Sigma$, from the category ${\mathcal
B}=\mathrm{Mod}$-${\bf B}$ to $\mathrm{Mod}$-${\bf T}$, is a ${\widetilde
{\mathbb C}}$-comodule functor, which is ${\widetilde {\mathbb C}}$-Galois
provided that
$\Sigma$ is a (not necessarily finite) Galois comodule.
The corresponding object in $\mathrm{RArr}(\mathrm{Mod}$-${\bf
A},\mathrm{Mod}$-${\bf B})$ is
$$
\big(\mathrm{Mod}\textrm{-}{\bf T},
((-)\ot_{\bf A}{\bf T},\mathrm{Hom}_{\bf T}({\bf T},-)),
((-)\ot_{\bf B}\Sigma,\mathrm{Hom}_{\bf T}(\Sigma,-)),
(-)\ot_{\bf A} {\bf C}, \beta^{-1}\big),
$$
where $\beta$ is the natural isomorphism, given in terms of the ${\bf
C}$-coaction $x \mapsto x_{(0)} \ot_{\bf A} x_{(1)}$ on $\Sigma$ as
$$
\mathrm{Hom}_{\bf T}(\Sigma,-) \ot_{\bf B} \Sigma \to (-) \ot_{\bf A} {\bf C},
\qquad
f \ot_{\bf A} x\mapsto f(x_{(0)}) \ot_{\bf A} x_{(1)}.
$$

In general, objects of the category $\overline{\mathrm{RArr}} ({\mathcal
  A},{\mathcal B})$ seem to have no interpretation similar to Corollary
\ref{cor:Gal.RArr}. Any object
$({\mathcal T},(N_A,R_A), (N_B,R_B), {\mathbb C},\xi)$ of
$\overline{\mathrm{RArr}} ({\mathcal A},{\mathcal B})$
determines a natural morphism
$$
\Phi:= (R_A \overline{\xi}^{-1}) \circ (\xi N_A) : C R_A N_A \to R_A N_A C,
$$
where $\overline{\xi}=(\epsilon^A N_B R_B N_A) \circ (N_A \xi N_A) \circ (N_A
C \eta^A)$ is an isomorphism by assumption. In terms of the structure
morphisms of the comonad ${\mathbb C}=(C,\Delta,\varepsilon)$, this morphism
$\Phi$ satisfies
\begin{eqnarray*}
&&(R_A N_A \Delta)\circ \Phi = (\Phi C) \circ (C \Phi) \circ (\Delta R_A N_A)
  \\
&&(R_A N_A \varepsilon)\circ \Phi = \varepsilon R_A N_A\\
&&\Phi \circ (C R_A \epsilon^A N_A) = (R_A \epsilon^A N_A C)\circ (R_A N_A
\Phi) \circ (\Phi R_A N_A)\\
&&\Phi \circ (C \eta^A ) =\eta^A C.
\end{eqnarray*}
Although these conditions are reminiscent to Definition \ref{def:entw} of a
mixed distributive law, in general $\Phi$ has no interpretation in terms of a
corresponding lifting
of ${\mathbb C}$ to the category of modules for the monad $R_A N_A$.
(But note that if both $\xi$ and ${\overline \xi}$
are isomorphisms then so is $\Phi$ and its inverse is a mixed distributive
law in the sense of Definition \ref{def:entw}.)

However, there are interesting (every-day-seen) examples of objects in
$\overline{\mathrm{RArr}}({\mathcal A},{\mathcal B})$ that can be interpreted
as Galois functors.
Namely, for two algebras ${\bf A}$ and ${\bf B}$, an ${\bf A}$-ring ${\bf T}$
and a ${\bf B}$-${\bf T}$ bimodule $\Sigma$, such that $\Sigma$ is a finitely
generated and projective right ${\bf T}$-module, there are adjunctions
\begin{eqnarray*}
&&\big( {\bf T}\ot_{\bf A} (-): {\bf A} \textrm{-}\mathrm{Mod} \to {\bf T}
\textrm{-}\mathrm{Mod},
{\bf T} \ot_{\bf T} (-): {\bf T} \textrm{-}\mathrm{Mod} \to {\bf A} \textrm{-}
\mathrm{Mod}\big)\\
&&\big( \Sigma^*\ot_{\bf B} (-): {\bf B} \textrm{-}\mathrm{Mod}\to {\bf T}
\textrm{-} \mathrm{Mod},
\Sigma \ot_{\bf T} (-) : {\bf T} \textrm{-} \mathrm{Mod}\to {\bf B} \textrm{-}
\mathrm{Mod} \big).
\end{eqnarray*}
Moreover, for any ${\bf A}$-coring ${\bf C}$,
$$
\big(\mathrm{Mod}\textrm{-} {\bf T},
((-)\ot_{\bf A} {\bf T},(-) \ot_{\bf T} {\bf T}),
((-)\ot_{\bf B} \Sigma,(-) \ot_{\bf T} \Sigma^*),
(-)\ot_{\bf A} {\bf C}, \xi = (-)\ot_{\bf T}\xi{\bf T}
\big)
$$
is an object of $\overline{\mathrm{RArr}} (\mathrm{Mod}\textrm{-} {\bf A},
\mathrm{Mod}\textrm{-} {\bf B})$ if and only if
$$
\big(
{\bf T} \textrm{-} \mathrm{Mod},
({\bf T}\ot_{\bf A} (-), {\bf T} \ot_{\bf T} (-)),
(\Sigma^*\ot_{\bf B} (-), \Sigma \ot_{\bf T} (-)),
{\bf C}\ot_{\bf A} (-),\overline {\xi} {\bf  A} \ot_{\bf T}(-)
\big)
$$
is an object of $\mathrm{RArr}({\bf A} \textrm{-}\mathrm{Mod} , {\bf B}
\textrm{-} \mathrm{Mod})$. That is, if and only if $\Sigma^* \ot_{\bf B}(-):
 {\bf B} \textrm{-}\mathrm{Mod}\to {\bf T} \textrm{-} \mathrm{Mod}$ is a
 Galois functor for a lifted comonad ${\bf C} \ot_{\bf A} (-)$ on ${\bf T}
 \textrm{-} \mathrm{Mod}$.

Thus, as objects of the category ${\mathrm{RArr}} ({\mathcal A},{\mathcal B})$
generalize (finite) right Galois comodules of corings arising from {\em right}
entwining structures, objects of $\overline{\mathrm{RArr}} ({\mathcal
  A},{\mathcal B})$ generalize, in a sense, finite left Galois comodules of
corings coming from {\em left} entwining structures.
Objects that belong to both categories
${\mathrm{RArr}} ({\mathcal A},{\mathcal B})$ and
$\overline{\mathrm{RArr}} ({\mathcal A},{\mathcal B})$
generalize dual pairs of finite left and right Galois comodules of the two
(isomorphic) corings arising from {\em bijective} entwining structures.

\subsection{Examples from monad morphisms}

More exotic examples of objects in a category $\mathrm{RArr}({\mathcal
  A},{\mathcal B})$ are obtained from {\em Galois extensions of monads
  by a comonad}.

Recall that in an adjunction $(L,R)$ -- with unit $\eta$
and counit $\epsilon$ --, the right adjoint $R$ is a right comodule functor
for some comonad ${\mathbb G}$, with coaction $g:R \to RG$ if and only if $L$
is a left ${\mathbb G}$-comodule functor via the coaction
\begin{equation}\label{eq:g.bar}
{\overline g}:=(\epsilon G L)\circ (L g L) \circ (L\eta): L\to G L.
\end{equation}

\begin{proposition}\label{prop:PREgalois}
Consider an adjoint pair of functors
$(L:{\mathcal K}\to {\mathcal C},R:{\mathcal C} \to {\mathcal K})$,
with unit $\eta$ and counit $\epsilon$ of the adjunction. Assume that $(R,g)$
is a right comodule functor for some comonad ${\mathbb G}$ on ${\mathcal C}$
(equivalently, $(L,{\overline g})$ is a left ${\mathbb G}$-comodule
functor, with coaction ${\overline g}$ in \eqref{eq:g.bar}). Then, if there
exists the equalizer
\begin{equation} \label{eq:B}
(B,\beta)=\mathrm{Equ}_{\mathrm{Fun}}(g L, R {\overline g}),
\end{equation}
then there exists a monad $\mathbb{B}=(B,m^{B},u^{B})$ on $\mathcal{K}$
such that $\beta$ gives rise to a morphism of monads
$\beta:\mathbb{B}\hookrightarrow RL$. Moreover, $\mathbb{B}$ is the unique
monad with this property.
\end{proposition}

\begin{proof}
By naturality and the adjunction relation, $(R \overline{g}) \circ \eta =
(gL)\circ \eta$. Thus by the universality of the equalizer \eqref{eq:B},
there exists a unique natural morphism $u^{B}:{\mathcal{K}}\rightarrow B$,
such that $\beta\circ u^{B}=\eta$.

Similarly, $(G\epsilon) \circ (\overline{g} R) = (\epsilon G)\circ (Lg)$. With
this identity at hand and using the fork property of \eqref{eq:B} (twice), one
checks that
$(R \epsilon L)\circ (\beta \beta)$ equalizes the parallel morphisms in
\eqref{eq:B}. Hence there exists a unique natural morphism $m^{B}:BB
\rightarrow  B$, such that $\beta\circ m^{B}=(R \epsilon L) \circ
(\beta\beta)$.
Associativity and unitality of the monad ${\mathbb B}:=(B,m^B,u^B)$ are obvious
by respective properties of the monad $(RL, R \epsilon L, \eta)$. The morphism
$\beta$ is compatible with the monad structures by construction.
\end{proof}

Consider a monad ${\bf T}$ on a category ${\mathcal K}$. We can apply
Proposition \ref{prop:PREgalois} to the associated adjunction $(L=F_{\mathbb
T},R= U_{\mathbb T})$, and a comonad ${\mathbb G}$ on ${\mathcal K}_{\mathbb
T}$. Note that the resulting notion of a ${\mathbb G}$-coaction $g$ on
$U_{\mathbb T}$ generalizes the notion of a {\em grouplike element} in a coring.
In this setting, if there exists the monad $\mathbb{B}=(B,m^{B},u^{B})$ on
$\mathcal{K}$ that is described in Proposition \ref{prop:PREgalois}, then it
will be denoted by $\mathbb{T}^{Co\left( G,g\right) }$ and it will be called
the \emph{coinvariant monad} of $\mathbb{T}$ with respect to $\mathbb{G}$ and
$g$. Note that the coinvariant monad $\mathbb{T}^{Co\left( G,g\right) }$
exists whenever ${\mathcal K}$ has equalizers.

\begin{proposition}\label{prop:galois}
Let $\mathcal{K}$ be a category with equalizers and coequalizers.
Let $\mathbb{T}=(T,m,u)$ be a monad on $\mathcal{K}$ and
$\mathbb{G=}\left( G,\Delta,\varepsilon\right) $ be a
comonad on $\mathcal{K}_{\mathbb{T}}$.
Assume that $T$ preserves coequalizers and that there exists a
right ${\mathbb G}$-coaction $g:U_{\mathbb T} \to U_{\mathbb T} G$.
Denote $\mathbb{B}:=\mathbb{T}^{Co\left( G,g\right) }$.
Let $(N_{B},R_{B})$ be the canonical adjunction, with unit $\eta^B$ and counit
$\epsilon^B$, associated as in Theorem \ref{theo:adjunct} to the canonical
inclusion  $\beta:\mathbb{ B\hookrightarrow T}$.
Under these assumptions, $N_B$
can be equipped with the structure of a left ${\mathbb G}$-comodule functor.
\end{proposition}

\begin{proof}
Consider the left ${\mathbb G}$-coaction
$\overline{g}=(\lambda^T G F_{\mathbb T})\circ (F_{\mathbb T} g F_{\mathbb T}) \circ
  (F_{\mathbb T} u^T):F_{\mathbb{T}}\rightarrow G F_{\mathbb{T}}$.
Recall that $F_{\mathbb T}$ is also a right ${\mathbb B}$-module functor, via
the action $f:=(\lambda^T F_{\mathbb T})\circ (F_{\mathbb T}\beta)$. Moreover,
the ${\mathbb B}$-action and the ${\mathbb G}$-coaction on $F_{\mathbb T}$
commute in the sense that $(Gf) \circ ({\overline g} B) = {\overline g}\circ
f$. This implies that $(G \chi^B)\circ (\overline{g} U_{\mathbb B})$
coequalizes the parallel morphisms (given by the ${\mathbb B}$-actions) in
$
(N_B, \chi^B)=\mathrm{Coequ}_{\mathrm{Fun}}\left(f
U_{\mathbb B}, F_{\mathbb T} U_{\mathbb B} \lambda^B \right).
$
So there exists a unique natural morphism
${\overline h}:N_{B}\rightarrow GN_{B}$, such that
\begin{equation}
{\overline h}\circ \chi ^{B}=(G \chi ^{B}) \circ (\overline{g}
{U_{\mathbb{B}}})\text{.}  \label{form:hGchiB}
\end{equation}
Coassociativity and counitality of ${\overline h}$ are obvious by the
analogous properties of $\overline g$. 
\end{proof}

The comodule functor ($N_B,{\overline h})$
in Proposition \ref{prop:galois} is a Galois functor provided that the
canonical comonad morphism
\begin{equation}
\mathrm{can}: = (G \epsilon ^{B}) \circ ({\overline h} {R_{B}}):
\left(N_{B}R_{B},N_{B}\eta^B {R_{B}},\epsilon ^{B}\right) \rightarrow
\left( G,\Delta ,\varepsilon \right)
\label{form:defcan}
\end{equation}
is an isomorphism.

\begin{theorem}
Let $\mathcal{K}$ be a category with equalizers and coequalizers
and let $\alpha :\mathbb{A}\to \mathbb{T}$
and $\beta:{\mathbb B} \to {\mathbb T}$
be morphisms of monads on $\mathcal{K}$.
Denote their canonical adjunctions (cf. Theorem \ref{theo:adjunct}) by
$(N_A,R_A)$ and $(N_B,R_B)$, respectively.
Assume that the underlying functors $A$ and $T$ preserve coequalizers.
Assume furthermore that the unit $\eta^B$ of the adjunction
$(N_B,R_B)$ is a regular natural monomorphism.
Then there is a bijective correspondence between the following sets of data.
\begin{itemize}
\item[{(a)}] objects
  $({\mathcal K}_{\mathbb T}, (N_A,R_A), (N_B,R_B), {\mathbb C}, \xi)$ of
  $\mathrm{RArr}({\mathcal K}_{\mathbb A}, {\mathcal K}_{\mathbb B})$,
\item[{(b)}] right coactions
$g:U_{\mathbb T} \to U_{\mathbb T} G$ for a lifting
of a comonad ${\mathbb C}$ on ${\mathcal K}_{\mathbb A}$
to a comonad ${\mathbb G}$ on ${\mathcal K}_{\mathbb T} \cong ({\mathcal
  K}_{\mathbb  A})_{R_AN_A}$, subject to the following conditions.
\begin{itemize}
\item[$\bullet$] ${\mathbb B}= {\mathbb T}^{Co(G,g)}$ and $\beta$ is the
  canonical inclusion ${\mathbb B}= {\mathbb T}^{Co(G,g)} \hookrightarrow
  {\mathbb T}$,
\item[$\bullet$] the canonical comonad morphism
\eqref{form:defcan} is an isomorphism.
\end{itemize}
\end{itemize}
\end{theorem}

\begin{proof}
By Theorem \ref{theo:main}, the module categories ${\mathcal K}_{\mathbb T}$
and $({\mathcal K}_{\mathbb A})_{R_A N_A}$ are isomorphic. Hence the
data in part (a) are in bijective correspondence with the objects
\begin{equation}\label{eq:prime}
(({\mathcal K}_{\mathbb A})_{R_A N_A}, (N'_A,R'_A), (N'_B,R'_B), {\mathbb C},
\xi')
\end{equation}
of $\mathrm{RArr}({\mathcal K}_{\mathbb A}, {\mathcal K}_{\mathbb
  B})$, where the primed functors are obtained by composing with the
isomorphism ${\mathcal K}_{\mathbb T}\cong ({\mathcal K}_{\mathbb A})_{R_A
N_A}$ on the appropriate side.
Note that $R'_A: ({\mathcal K}_{\mathbb A})_{R_A N_A}\to {\mathcal K}_{\mathbb
  A}$ is the forgetful functor corresponding to the monad $R_A N_A$.
Therefore, by Theorem \ref{thm:Gal.com.tr} there is a bijection between the
data in \eqref{eq:prime} and
liftings of ${\mathbb C}$ to a comonad ${\widetilde{\mathbb C}}$ on
$({\mathcal K}_{\mathbb A})_{R_A N_A}$ together with a ${\widetilde{\mathbb
C}}$-Galois structure on the functor $N'_B$.
The isomorphism ${\mathcal K}_{\mathbb T}\cong ({\mathcal K}_{\mathbb A})_{R_A
N_A}$ takes a comonad ${\widetilde{\mathbb C}}$ on $({\mathcal K}_{\mathbb
A})_{R_A N_A}$ to a comonad ${\mathbb G}$ on ${\mathcal K}_{\mathbb
T}$. Clearly, $N_B$ is a ${\widetilde{\mathbb C}}$-Galois functor if and only
if $N'_B$ is a ${\mathbb G}$-Galois functor.
Thus the data in part (a) are in bijective correspondence with liftings of
${\mathbb C}$ to a comonad ${\mathbb G}$ on ${\mathcal K}_{\mathbb T}$
together with a ${\mathbb G}$-Galois structure on $N_B$.

The data in part (b) determine a ${{\mathbb G}}$-Galois structure on $N_B$ by
Proposition \ref{prop:galois}.
Conversely, a ${{\mathbb G}}$-Galois
functor $(N_B,{\overline h})$ determines a right ${\mathbb G}$-coaction
$h:=(R_B G \epsilon^B)\circ (R_B {\overline h} R_B)\circ (\eta^B R_B)$ on
$R_B$ and a right ${\mathbb G}$-coaction
$g:= U_{\mathbb B} h=(U_{\mathbb T} \mathrm{can}) \circ (U_{\mathbb B} \eta^B R_B)$
on $U_{\mathbb T}=U_{\mathbb B} R_B$,
where $\mathrm{can}$ is the isomorphism \eqref{form:defcan}.
It remains to show that for this coaction $g$,
\begin{eqnarray}
(B,\beta )&=&\mathrm{Equ}_{\mathrm{Fun}}(g {F_{\mathbb{T}}},(U_{\mathbb{T}}
  \lambda^T
{GF_{\mathbb{T}}}) \circ  (T g {F_{\mathbb{T}}})\circ (T u^{T}))\nonumber\\
&\equiv&\mathrm{Equ}_{\mathrm{Fun}}((U_{\mathbb{B}}\eta^B
{R_{B}F_{\mathbb{T}}}),(U_{\mathbb{T}}
\lambda^T {N_{B}R_{B}F_{\mathbb{T}}}) \circ (TU_{\mathbb B} \eta^B
       {R_{B}F_{\mathbb{T}}})\circ (Tu^{T})) \label{eq:B_eq}.
\end{eqnarray}
Since both $F_{\mathbb T}$ and $N_B F_{\mathbb B}$ are left adjoints of
$U_{\mathbb T}= U_{\mathbb B} R_B$, there is a natural isomorphism $\gamma:
N_B F_{\mathbb B} \to F_{\mathbb T}$. Clearly, \eqref{eq:B_eq} is equivalent
to
\begin{eqnarray}\label{eq:beta_eq}
&&\left(B, (U_{\mathbb T} \gamma^{-1})\circ \beta \right) = \nonumber \\
&&\mathrm{Equ}_{\mathrm{Fun}}\left(
(U_{\mathbb T} N_B R_B \gamma^{-1}) \circ (U_{\mathbb B} \eta^B R_B F_{\mathbb
  T}) \circ (U_{\mathbb T} \gamma) \right., \nonumber\\
&&\qquad \qquad \left. (U_{\mathbb T} N_B R_B \gamma^{-1}) \circ
  (U_{\mathbb{T}} \lambda^T {N_{B}R_{B}F_{\mathbb{T}}}) \circ (TU_{\mathbb B}
  \eta^B {R_{B}F_{\mathbb{T}}})\circ (Tu^{T}) \circ (U_{\mathbb T} \gamma)
  \right).
\end{eqnarray}
Using the explicit form of
$\gamma =(\epsilon^B F_{\mathbb{T}})\circ (N_{B} \lambda^B
R_{B}F_{\mathbb{T}}) \circ (N_{B} F_{\mathbb{B}} u^{T})$, it is
straightforward to check that
\begin{eqnarray*}
&&(U_{\mathbb T} \gamma^{-1})\circ \beta = U_{\mathbb B} \eta^B F_{\mathbb
    B},\\
&&(U_{\mathbb T} N_B R_B \gamma^{-1}) \circ (U_{\mathbb B} \eta^B R_B F_{\mathbb
  T}) \circ (U_{\mathbb T} \gamma)= U_{\mathbb B} \eta^B R_B N_B F_{\mathbb
    B},\\
&&(U_{\mathbb T} N_B R_B \gamma^{-1}) \circ
  (U_{\mathbb{T}} \lambda^T {N_{B}R_{B}F_{\mathbb{T}}}) \circ (TU_{\mathbb B}
  \eta^B {R_{B}F_{\mathbb{T}}})\circ (Tu^{T}) \circ (U_{\mathbb T} \gamma) =
U_{\mathbb B}  R_B N_B \eta^B F_{\mathbb B}.
\end{eqnarray*}
Thus \eqref{eq:beta_eq} is equivalent to
\begin{equation}\label{eq:gamma_f}
\left(B, U_{\mathbb B} \eta^B F_{\mathbb B}\right)=
\mathrm{Equ}_{\mathrm{Fun}}\left(U_{\mathbb B} \eta^B R_B N_B F_{\mathbb
  B},U_{\mathbb B}  R_B N_B \eta^B
F_{\mathbb B}\right).
\end{equation}
By the assumption that $\eta^B$ is an equalizer,
Lemma \ref{lem:eta.reg} and Lemma \ref{Lem: Funct1}, we conclude that
$\eta^B F_{\mathbb B}$ is the equalizer of
the parallel morphisms $\eta^B R_B N_B F_{\mathbb B}$ and $ R_B N_B \eta^B
F_{\mathbb B}$. Since $U_{\mathbb B}$ is a right adjoint functor, it preserves
equalizers. This proves \eqref{eq:gamma_f} hence \eqref{eq:B_eq}.

Bijectivity of the constructed maps between the data in part (a) and part (b)
is checked by a straightforward computation,
making use of \eqref{form: epsil} and \eqref{form:UAdietaA} .
\end{proof}

\section{Equivalence between regular comonad arrows and pre-torsors}
\label{sec:pre-tor}

In this section we study the category $\mathrm{RArr}({\mathcal A},{\mathcal B})$
introduced in Section \ref{sec:entw}
(possessing the subcategory $\mathrm{Gal}({\mathcal A},{\mathcal B})$). Our
main aim is to find a full subcategory $\mathrm{RArr}_{reg}({\mathcal
A},{\mathcal B})$ that is equivalent to a full subcategory of the category of
pre-torsors, introduced below.
The subcategory $\mathrm{RArr}_{reg}({\mathcal A},{\mathcal B})$ has some
intersection with the subcategory $\mathrm{Gal}({\mathcal
A},{\mathcal B})$ of $\mathrm{RArr}({\mathcal A},{\mathcal B})$. Hence the
intersection (let's denote it by $\mathrm{Gal}_{reg} ({\mathcal A},{\mathcal
  B})$), is equivalent to a subcategory of the category of pre-torsors.
Since any torsor corresponding to a faithfully flat Hopf Galois object
belongs to  $\mathrm{Gal}_{reg} ({\mathcal A},{\mathcal B})$,
the results of this section not only generalize the relation between
faithfully flat Hopf Galois objects and torsors over commutative rings, but
also yield a deeper explanation of it.
It is an open question how to interpret those objects of
$\mathrm{RArr}_{reg}({\mathcal A},{\mathcal B})$ that do not belong to
$\mathrm{Gal}({\mathcal A},{\mathcal B})$.

\subsection{Pre-torsors}

In this section we introduce a further key notion of these notes, {\em
  pre-torsors} over two adjunctions.
They provide examples of herd functors in \cite[Appendix]{BV}.

\begin{definition} \label{def:pre-tors}
{\em Objects} in the category $\mathrm{PreTor}({\mathcal A}, {\mathcal B})$,
called {\em pre-torsors}, consist of an object $({\mathcal
  T},(N_A,R_A),(N_B,R_B))$ of
$\mathrm{Adj}({\mathcal A}, {\mathcal B})$ together with a natural morphism
$\tau :R_{A}N_{B}\rightarrow R_{A}N_{B}R_{B}N_{A}R_{A}N_{B}$, subject to the
following conditions.
\begin{itemize}
\item $(R_{A}N_{B}R_{B} \epsilon^A {N_{B}}) \circ \tau
=R_{A}N_{B}\eta ^{B} $,
\item $(R_{A}\epsilon^B {N_{A}R_{A}N_{B}}) \circ \tau =\eta^A {R_{A}N_{B}}$,
\item $(R_{A}N_{B}R_{B}N_{A} \tau) \circ \tau =(\tau
{R_{B}N_{A}R_{A}N_{B}})\circ \tau ,$
\end{itemize}
where $\eta^A$ and $\epsilon^A$ denote the unit and counit of the
adjunction $(N_A,R_A)$ and
analogous notations $\eta^B$ and $\epsilon^B$ are used for
$(N_B,R_B)$.

A {\em morphism} of pre-torsors
$$
({\mathcal T},(N_A,R_A),(N_B,R_B),\tau) \to
({\mathcal T}',(N'_A,R'_A),(N'_B,R'_B),\tau')
$$
is a morphism $F:({\mathcal T},(N_A,R_A),(N_B,R_B)) \to ({\mathcal
  T}',(N'_A,R'_A),(N'_B,R'_B))$ in $\mathrm{Adj}({\mathcal A},{\mathcal B})$,
such that
\begin{equation}\label{eq:tor_mor}
(R_A b R_B a R_A b) \circ \tau = \tau' \circ (R_A b),
\end{equation}
where the natural morphisms $a$ and $b$ were defined in \eqref{eq:a&b}.
\end{definition}

Note that $R_A b$ is a natural morphism $R_A N_B \to R'_A N'_B$ and $R_B a$
is a natural morphism $R_B N_A \to R'_B N'_A$.
We leave it to the reader to check that the composite of two morphisms
of pre-torsors is a morphism of pre-torsors again.

\begin{example}\label{ex:tor_bim}
For a commutative ring $k$, consider $k$-algebra homomorphisms $\boldsymbol{
  \alpha:}$ $\boldsymbol{A\to T}$ and $\boldsymbol{\beta :B\to T}$. Denote by
$(N_A,R_A)$
(resp. $(N_B,R_B)$) the corresponding `extension of scalars' and `restriction
of scalars' functors between the module categories of these algebras.
A pre-torsor $\boldsymbol{(\alpha,\beta,\tau)}$,
over $k$-algebras $\boldsymbol{A}$ and $\boldsymbol{B}$ as in \cite[Definition
  3.1]{BT}, induces a pre-torsor
$(\mathrm{Mod}\textrm{-}{\boldsymbol T},(N_A,R_A), (N_B,R_B),  (-)
\otimes_{\boldsymbol{B}} \boldsymbol{\tau})$ in the sense of Definition
\ref{def:pre-tors}.
\end{example}

\begin{example}\label{ex:herd}
More generally, for two associative and unital algebras ${\bf A}$ and ${\bf
  B}$, consider
an ${\bf A}$-ring ${\bf T}$ and a ${\bf B}$-${\bf T}$ bimodule $\Sigma$. As
we have seen in Section \ref{sec:bimod.ex}, these data determine two pairs of
adjoint functors:
\begin{eqnarray*}
&N_A = (-)\ot_{\bf A} {\bf T}: \mathrm{Mod}\textrm{-}{\bf A} \to
  \mathrm{Mod}\textrm{-}{\bf T}
\qquad
R_A=
\mathrm{Hom}_{\bf T}({\bf T},-): \mathrm{Mod}\textrm{-}{\bf T} \to
\mathrm{Mod}\textrm{-}{\bf A}\\
&\ N_B = (-)\ot_{\bf B} \Sigma : \mathrm{Mod}\textrm{-}{\bf B} \to
  \mathrm{Mod}\textrm{-}{\bf T}
\qquad
R_B=
\mathrm{Hom}_{\bf T}(\Sigma,-): \mathrm{Mod}\textrm{-}{\bf T} \to
\mathrm{Mod}\textrm{-}{\bf B}.
\end{eqnarray*}
Then $R_A N_B= (-)\ot_{\bf B}\Sigma$. If $\Sigma$ is finitely generated and
projective as a right ${\bf T}$-module, then $R_B N_A = (-) \ot_A \Sigma^*$,
where the notation $\Sigma^*:=\mathrm{Hom}_{\bf T} (\Sigma,{\bf T})$ is
used. Hence a pre-torsor is induced by any ${\bf B}$-${\bf A}$-bimodule map
$$
{\boldsymbol{\tau}}: \Sigma \to \Sigma\ot_{\bf A} \Sigma^*\ot_{\bf B} \Sigma,
\qquad
x\mapsto x^{\langle 1\rangle} \ot_{\bf A} x_{\langle 2\rangle} \ot_{\bf B}
x^{\langle 3\rangle} ,
$$
(with implicit summation understood), which is subject to the following
conditions, for all $x\in \Sigma$.
\begin{eqnarray*}
&x^{\langle 1\rangle} x_{\langle 2\rangle} (-) \ot_{\bf B} x^{\langle
3\rangle} = \mathrm{Id}_\Sigma\ot_{\bf B} x &\in \mathrm{End}_{\bf
T}(\Sigma) \ot_{\bf B} \Sigma,\\
&x^{\langle 1\rangle} \ot_{\bf A} x_{\langle 2\rangle}(x^{\langle 3\rangle}) =
x \ot_{\bf A} 1_{\bf T} &\in \Sigma \ot_{\bf A} {\bf T},\\
&x^{\langle 1\rangle} \ot
x_{\langle 2\rangle}\ot
x^{\langle 3\rangle \langle 1\rangle} \ot
{x^{\langle 3\rangle}}_{\langle 2\rangle} \ot
x^{\langle 3\rangle \langle 3\rangle} &=
x^{\langle 1\rangle \langle 1\rangle} \ot
{x^{\langle 1\rangle}}_{\langle 2\rangle} \ot
x^{\langle 1\rangle \langle 3\rangle} \ot
x_{\langle 2\rangle} \ot
x^{\langle 3\rangle} \\
& \qquad\qquad \qquad &\in \Sigma\ot_{\bf A} \Sigma^*\ot_{\bf B}
\Sigma  \ot_{\bf A} \Sigma^*\ot_{\bf B} \Sigma.
\end{eqnarray*}
This structure is an example of a {\em bimodule} herd, introduced in
\cite[Definition 2.4]{BV}, over the ring maps ${\bf A}\to {\bf T}$ and ${\bf
  B} \to \mathrm{End}_{\bf T}(\Sigma)$.
\end{example}

\subsection{From regular comonad arrows to pre-torsors} \label{sec:Gal>Tor}

Our next aim is to find an equivalence between (certain) objects in
$\mathrm{RArr}({\mathcal A}, {\mathcal B})$ on one hand, and (certain)
pre-torsors on the other hand. In the same way as it happens with the relation
between Galois extensions by corings and pre-torsors over rings, we will see
in this section that any object in $\mathrm{RArr}({\mathcal A},{\mathcal B})$
(or $\overline{\mathrm{RArr}}({\mathcal B},{\mathcal A})$) determines a
pre-torsor.
In the next section we ask about a converse construction and look for
conditions under which a pre-torsor
determines an object in $\mathrm{RArr}({\mathcal A}, {\mathcal B})$ (or
$\overline{\mathrm{RArr}}({\mathcal B},{\mathcal A})$). Our
main result is a category equivalence in Section \ref{sec:cat_eq}.

\begin{theorem} \label{thm:Arr>Tor}
(1) Consider an object $({\mathcal T},(N_A,R_A),(N_B,R_B), {\mathbb C},\xi)$
  in $\mathrm{RArr}({\mathcal A},{\mathcal B})$. Then
  $({\mathcal T},(N_A,R_A),(N_B,R_B),\tau)$ is a pre-torsor, with pre-torsor
  map
\begin{equation*}
\tau =(\xi N_A R_A N_B)\circ (C \eta^A R_A N_B) \circ (\xi^{-1} N_B) \circ
(R_A N_B \eta^B) .
\end{equation*}
(2) For a morphism
$$
(F,t):
({\mathcal T},(N_A,R_A),(N_B,R_B), {\mathbb C},\xi) \to
({\mathcal T}',(N'_A,R'_A),(N'_B,R'_B), {\mathbb C}',\xi')
$$
in $\mathrm{RArr}({\mathcal A},{\mathcal B})$, $F$ is a morphism of
  the corresponding pre-torsors in part (1).
\end{theorem}

\begin{proof}
Part (1) is verified by a direct computation, using the definition of an
object in $\mathrm{RArr}({\mathcal A},{\mathcal B})$, naturality and the
adjunction relations. Similarly, to check claim (2), one has to use
the definition of a morphism in $\mathrm{RArr}({\mathcal A},{\mathcal B})$,
together with naturality and the adjunction relations.
\end{proof}

Symmetrically to Theorem \ref{thm:Arr>Tor}, we have
\begin{theorem} \label{thm:LArr>Tor}
(1) Consider an object $({\mathcal T},(N_B,R_B),(N_A,R_A), {\mathbb D},\zeta)$
  in $\overline{\mathrm{RArr}}({\mathcal B},{\mathcal A})$. Then
  $({\mathcal T},(N_A,R_A),(N_B,R_B),\tau)$ is a pre-torsor, with pre-torsor
  map
\begin{equation*}
\tau = (R_A N_B R_B {\overline \zeta})\circ (R_A N_B \eta^B D) \circ (R_A
{\overline  \zeta}^{-1}) \circ (\eta^A R_A N_B) ,
\end{equation*}
where ${\overline \zeta}:= (\epsilon^B N_A R_A N_B) \circ (N_B \zeta N_B) \circ
(N_B D \eta^B)$.

(2) For a morphism
$$
(F,t):
({\mathcal T},(N_B,R_B),(N_A,R_A), {\mathbb D},\zeta) \to
({\mathcal T}',(N'_B,R'_B),(N'_A,R'_A), {\mathbb D}',\zeta')
$$
in $\overline{\mathrm{RArr}}({\mathcal B},{\mathcal A})$, $F$ is a morphism of
  the corresponding pre-torsors in part (1).
\end{theorem}

In Theorem \ref{thm:Arr>Tor} we constructed, in fact, a functor from
$\mathrm{RArr}({\mathcal A},{\mathcal B})$ to $\mathrm{PreTor}({\mathcal
  A},{\mathcal B})$ and in Theorem \ref{thm:LArr>Tor} we constructed a functor
from $\overline{\mathrm{RArr}}({\mathcal B},{\mathcal A})$ to
$\mathrm{PreTor}({\mathcal
  A},{\mathcal B})$. In later sections of these notes we want to see on which
objects are these functors equivalences. The problem will be
divided to two steps. In Section
\ref{sec:Tor>Gal}, we investigate on what subcategory of the
category of pre-torsors we can define functors to
$\mathrm{RArr}({\mathcal A},{\mathcal B})$ and
$\overline{\mathrm{RArr}}({\mathcal B},{\mathcal A})$.
After that, in Section
\ref{sec:cat_eq}, we prove that the functors in Sections
\ref{sec:Gal>Tor} and \ref{sec:Tor>Gal} give rise to inverse equivalences
between appropriate subcategories, indeed.

\subsection{From pre-torsors to regular comonad arrows} \label{sec:Tor>Gal}

The aim of this section is to find criteria on a pre-torsor, under
which it determines an object in $\mathrm{RArr}({\mathcal
  A}, {\mathcal B})$ (or $\overline{\mathrm{RArr}}({\mathcal B},{\mathcal
  A})$). The main result of the section is the following.

\begin{theorem}\label{thm:Tor>Arr}
Consider two categories ${\mathcal A}$ and ${\mathcal B}$ both of which
possess equalizers, and an object $({\mathcal T},(N_A,R_A),(N_B,R_B),\tau)$ in
$\mathrm{PreTor}({\mathcal A},{\mathcal B})$, such that
\begin{itemize}
\item the unit $\eta^A$ of the adjunction $(N_A,R_A)$ and the unit $\eta^B$ of
the adjunction $(N_B,R_B)$ are regular natural monomorphisms;
\item the functors $N_A$ and $N_B$ preserve equalizers.
\end{itemize}
Under these assumptions, the following assertions hold.

(1) The equalizer $ \left( C,i\right)
  =\mathrm{Equ}_{\mathrm{Fun}}(\omega^{l},\omega^{r})$ of the natural
  morphisms
$$
\omega ^{l}:=(R_{A}N_{B}R_{B}N_{A}R_{A} \epsilon^B {N_{A}})
\circ (\tau {R_{B}N_{A}})\qquad \textrm{and}\qquad
\omega ^{r}:=R_{A}N_{B}R_{B}N_{A} \eta ^{A}
$$
defines a comonad ${\mathbb C}=(C,\Delta^C,\varepsilon^C)$ on ${\mathcal A}$
such that the functor $C$ preserves equalizers.

(2) There is an object
$({\mathcal T},(N_A,R_A),(N_B,R_B),{\mathbb C},\xi)\in \mathrm{RArr}({\mathcal
  A},{\mathcal B})$, where the comonad ${\mathbb C}$ was constructed in part
  (1).

(3) For a morphism $F$, between pre-torsors
  $({\mathcal T},(N_A,R_A),(N_B,R_B),\tau)$ and $({\mathcal
  T}',$ $(N'_A,R'_A),(N'_B,R'_B),\tau')$
  both of which satisfy the conditions in the theorem, there exists a unique
  morphism $t:{\mathbb C} \to {\mathbb C}'$ between
       the comonads in part (1) such that $(F,t)$ is a morphism
       between the objects of $\mathrm{RArr}({\mathcal
  A},{\mathcal B})$ in part (2).
\end{theorem}

Before turning to prove Theorem \ref{thm:Tor>Arr}, we give a motivating
example of a situation in which the assumptions of the theorem hold.

\begin{example}
Consider an object
$$
(\mathrm{Mod}\textrm{-}{\boldsymbol T}, (-\ot_A {\boldsymbol T},
\mathrm{Hom}_{\bf T}({\bf T},-)),
(-\ot_B {\boldsymbol T}, \mathrm{Hom}_{\bf T}({\bf T},-)),
- \otimes_{\boldsymbol{B}}
\boldsymbol{\tau})
$$
in $\mathrm{PreTor}(\mathrm{Mod}\textrm{-}{\boldsymbol A},
\mathrm{Mod}\textrm{-}{\boldsymbol B})$
as in Example \ref{ex:tor_bim}, induced by a pre-torsor $(\boldsymbol{
  \alpha,\beta,\tau})$ in \cite[Definition
  3.1]{BT}. If ${\boldsymbol T}$ is a faithfully flat left
$\boldsymbol{A}$-module (via $\boldsymbol{\alpha}$) and a faithfully flat left
$\boldsymbol{B}$-module (via $\boldsymbol{\beta}$) then all assumptions in
Theorem \ref{thm:Tor>Arr} hold.

In the more general situation in Example \ref{ex:herd}, the assumptions in
Theorem \ref{thm:Tor>Arr} hold provided that ${\bf T}$ is a faithfully flat
left ${\bf A}$-module and $\Sigma$ is a faithfully flat left ${\bf
  B}$-module.
\end{example}

Part (1) of Theorem \ref{thm:Tor>Arr} holds in a more general situation in the
following lemma. 

\begin{lemma}\label{lem:C&D}
Consider two categories ${\mathcal A}$ and ${\mathcal B}$ possessing
equalizers, and equalizer preserving functors as in the (non-commutative)
diagram
$$
\xymatrix{
{\mathcal A}\ar@(ul,ur)[rr]^-{P}\ar@(dl,ul)^-{R}&&
{\mathcal B}\ar@(dr,dl)[ll]^-{Q}\ar@(dr,ur)_-{S}
} .
$$
Let
$$
r:{\mathcal A}\to R,\qquad
s:{\mathcal B} \to S,\qquad
w: QP \to R,\qquad
z: PQ \to S,\qquad
\tau: Q \to QPQ
$$
be natural morphisms, subject to the following conditions.
\begin{itemize}
\item[{(i)}] $r = \mathrm{Equ}_{\mathrm{Fun}}(rR,Rr)$ and
$s = \mathrm{Equ}_{\mathrm{Fun}}(sS,Ss)$;
\item[{(ii)}] $(QP\tau)\circ \tau = (\tau PQ) \circ \tau$;
\item[{(iii)}] $(Qz)\circ \tau=Qs$ and $(wQ)\circ \tau= rQ$.
\end{itemize}
Then the following assertions hold.

(1) There is a comonad ${\mathbb C}=(C,\Delta^C,\varepsilon^C)$ on ${\mathcal
  A}$, such that $C$ preserves equalizers and $Q$ is a left ${\mathbb
  C}$-comodule functor;

(2) There is a comonad ${\mathbb D}=(D,\Delta^D,\varepsilon^D)$ on ${\mathcal
  B}$, such that $D$ preserves equalizers and $Q$ is a right ${\mathbb
  D}$-comodule functor;

(3) $Q$ is a ${\mathbb C}$-${\mathbb D}$ bicomodule functor.
\end{lemma}

\begin{proof}
(1)
Consider the equalizer of natural morphisms
\begin{equation}\label{eq:C.eq}
\xymatrix{
C \ar[r]^-{i}&
QP \ar@<2pt>[r]^-{\omega^l} \ar@<-2pt>[r]_-{\omega^r}&
QPR,
}
\end{equation}
where
$\omega^l=(QPw)\circ (\tau P)$ and $\omega^r=QPr$.
(This equalizer exists by Lemma \ref{Lem: Equ}.)
By Corollary \ref{cor:pres_eq}, $C$ preserves equalizers.
By assumptions (ii) and (iii), $(\omega^l Q)\circ \tau = (\omega^r Q)\circ
\tau$. Hence by Lemma \ref{Lem: Funct1},
there is a unique natural morphism $c: Q \to CQ$ such that
\begin{equation}\label{eq:C-coac}
(iQ)\circ c =\tau.
\end{equation}
Moreover, use \eqref{eq:C-coac}, assumption (ii), definition of the morphism
$i$ and naturality to derive
$$
(iQPR)\circ (C\omega^l)\circ (cP)\circ i =
(iQPR)\circ (C\omega^r)\circ (cP)\circ i .
$$
Since $iQPR$ is a monomorphism by Lemma \ref{Lem: Funct1},
and $Ci$ is the equalizer of $C\omega^l$ and
$C\omega^r$, this implies the
existence of a unique natural morphism $\Delta^C:C \to CC$, such that
$(Ci)\circ \Delta^C = (cP)\circ i$. Equivalently,
\begin{equation}\label{form:DeltaCbar}
(ii) \circ \Delta^C = (\tau P)\circ i.
\end{equation}
Finally, by the definition of the morphism $i$, assumption (iii) and
naturality, 
$$
(Rr)\circ w\circ i =(rR)\circ w\circ i.
$$
By assumption (i) this implies the existence of a unique natural morphism
$\varepsilon ^C:C\to {\mathcal A}$ such that
\begin{equation}\label{form:epsilonC}
r\circ \varepsilon^C = w\circ i.
\end{equation}
By \eqref{form:DeltaCbar} and assumption (ii), $\Delta^C$ is coassociative. By
\eqref{form:DeltaCbar}, \eqref{form:epsilonC} and assumption (iii) on one
hand, and definition of $i$ via \eqref{eq:C.eq} on the other hand,
$\varepsilon^C$ is counit of $\Delta^C$. By \eqref{eq:C-coac},
\eqref{form:DeltaCbar} and assumption (ii), $c$ is a coassociative coaction.
By \eqref{eq:C-coac}, \eqref{form:epsilonC} and assumption (iii), it is
counital.

Part (2) is proven symmetrically. That is, we define an (equalizer preserving)
functor $D$ as the equalizer
\begin{equation}\label{eq:D.eq}
\xymatrix{
D \ar[r]^-{j}& PQ \ar@<2pt>[r]^-{\theta^l} \ar@<-2pt>[r]_-{\theta^r}&
SPQ,
}
\end{equation}
where $\theta^l=(zPQ)\circ (P\tau)$ and $\theta^r=sPQ$. The to-be-coaction
$d:Q \to QD$ is the unique morphism such that
\begin{equation}\label{eq:D.coac}
(Qj)\circ d =\tau.
\end{equation}
The coproduct and counit of the comonad ${\mathbb D}$ are defined via the
respective conditions
\begin{equation}\label{eq:D.str}
(jj)\circ \Delta^D = (P\tau)\circ j \qquad \textrm{and}\qquad s\circ
\varepsilon^D = z\circ j.
\end{equation}

(3) The identity $(Cd)\circ c = (c D)\circ d$ follows by \eqref{eq:C-coac},
\eqref{eq:D.coac} and assumption (ii).
\end{proof}

\begin{lemma} \label{lem:xi}
Under the assumptions and using the notations of Theorem \ref{thm:Tor>Arr},
$$
\xi:=(R_{A}N_{B}R_{B} \epsilon ^{A})
\circ (i {R_{A}}):C R_{A}\rightarrow R_{A}N_{B}R_{B}
$$
is a natural isomorphism.
\end{lemma}

\begin{proof}
It follows easily by the pre-torsor axioms that
\begin{equation}
(\omega^l {R_{A}})\circ (R_{A}N_{B}R_{B}N_{A}R_{A} \epsilon ^{B}) \circ (\tau
  {R_{B}})=(\omega^r {R_{A}})\circ
(R_{A}N_{B}R_{B}N_{A}R_{A} \epsilon ^{B}) \circ (\tau
  {R_{B}}).  \label{form:necessaria}
\end{equation}
Since the category ${\mathcal A}$ has equalizers by assumption, it follows by
Lemma \ref{Lem: Funct1} that $i R_A$ is the equalizer of $\omega^l
R_A$ and $\omega^r R_A$.
So by its universality there exists a unique natural morphism
$\xi^{\prime}: R_{A}N_{B}R_{B}\rightarrow C R_{A}$, such that
\begin{equation}
(i {R_{A}})\circ \xi^{\prime }=(R_{A}N_{B}R_{B}N_{A}R_{A}\epsilon
  ^{B}) \circ (\tau {R_{B}}).
  \label{form:phiprimo}
\end{equation}
A direct computation verifies that $\xi^{\prime }$ is the two-sided
inverse of $\xi$.
\end{proof}

We are ready to prove Theorem \ref{thm:Tor>Arr}, which is the main result of
the section.

\begin{proof}[Proof of Theorem \ref{thm:Tor>Arr}]
(1) The functors $R_A$ and $R_B$ are right adjoints, hence they
  preserve equalizers. The left adjoint functors $N_A$ and $N_B$ preserve
  equalizers by assumption. Hence
part (1) follows by substituting in Lemma \ref{lem:C&D} $Q=R_A N_B$, $P= R_B
N_A$, $R=R_A N_A$, $S= R_B N_B$, $r=\eta^A$, $s=\eta^B$, $w= R_A \epsilon^B
N_A$ and $z= R_B \epsilon^A N_B$.

(2) We claim that $({\mathcal T},(N_A,R_A),(N_B,R_B),{\mathbb C},\xi)$ is an
object in $\mathrm{RArr}({\mathcal A},{\mathcal B})$, where ${\mathbb C}$ is
the comonad in part (1) and $\xi$ is the isomorphism
constructed in Lemma \ref{lem:xi}. For that, we need to show that $\xi$
satisfies conditions \eqref{eq:com.tr}. This is a consequence of its
construction, naturality, \eqref{form:DeltaCbar} and \eqref{form:epsilonC},
the pre-torsor axioms and the adjunction relations.

(3) Take a morphism
$$
F:({\mathcal T}, (N_A,R_A), (N_B,R_B),\tau) \to
({\mathcal T}', (N'_A,R'_A), (N'_B,R'_B),\tau')
$$
of pre-torsors, with corresponding natural morphisms $a$ and $b$ in
\eqref{eq:a&b}.
For the morphisms $\omega^l$ and $\omega^r$, defined for the pre-torsor
$({\mathcal T}, (N_A,R_A),$ $ (N_B,R_B),\tau)$
in part (1), and the analogous morphisms $
\omega^{\prime l}$ and
$\omega^{\prime r}$ for $({\mathcal T}', (N'_A,R'_A),$ $ (N'_B,R'_B),\tau')$, we
claim that
\begin{eqnarray}\label{eq:omega_comp}
&& \omega^{\prime l} \circ (R_A b R_B a) = (R_A b R_B a R_A a) \circ
  \omega^l\qquad \textrm{and}\\
&& \omega^{\prime r} \circ (R_A b R_B a) = (R_A b R_B a R_A a) \circ
  \omega^r.\nonumber
\end{eqnarray}
The first identity follows by the explicit form of the morphisms $\omega^l$ and
$\omega^{\prime l}$, combining with the fact that $F$ is a morphism of
pre-torsors. In light of the form of the morphisms $\omega^r$ and
$\omega^{\prime r}$, the second identity is a
consequence of naturality and the compatibility of $a$ with the units of the
adjunctions $(N_A,R_A)$ and $(N'_A,R'_A)$, cf. \eqref{eq:a&b}. Composing both
identities in \eqref{eq:omega_comp} by the equalizer
$(C,i):=\mathrm{Equ}_{\mathrm{Fun}}(\omega^l, \omega^r)$, we obtain
$$
 \omega^{\prime r} \circ (R_A b R_B a) \circ i =
\omega^{\prime l} \circ (R_A b R_B a) \circ i.
$$
Thus by universality of the equalizer
$(C',i'):=\mathrm{Equ}_{\mathrm{Fun}}(\omega^{\prime l}, \omega^{\prime r})$,
there exists a natural morphism $t:C\rightarrow C^{\prime }$ such
that
\begin{equation}\label{form:jutile}
i^{\prime }\circ t=(R_A b R_B a)\circ i\text{.}
\end{equation}
We prove next that $t$ is a comonad morphism.
Using definitions \eqref{form:epsilonC} of the counit $\varepsilon^C$ and
\eqref{form:jutile} of the morphism $t$ together with \eqref{eq:a&b},
one checks that $\eta ^{^{\prime
  }A}\circ \varepsilon ^{C}=\eta
^{^{\prime }A}\circ \varepsilon ^{C^{\prime }}\circ t$.
So by the monomorphism property of $\eta ^{^{\prime }A}$, $\varepsilon
^{C}=\varepsilon ^{C^{\prime }}\circ t$.
Using definitions \eqref{form:jutile} of $t$
and \eqref{form:DeltaCbar} of the coproduct $\Delta^C$ together with
the fact that $F$ is a morphism of pre-torsors,
we deduce that
$$
(i^{\prime }i^{\prime })\circ (t t)\circ \Delta ^{C} =
(i^{\prime } i^{\prime })\circ \Delta^{C^{\prime }}\circ t\text{.}
$$
Since $i'i'= (R_A' N_B' R_B' N_A' i')\circ (i' C')$ is monic, we conclude that
$(t t)\circ \Delta ^{C}=\Delta ^{C^{\prime }}\circ t$. Condition
\eqref{eq:RGal_mor} follows by constructions of $\xi'$ and $t$, and
\eqref{eq:a&b}.
\end{proof}

Symmetrically to Theorem \ref{thm:Tor>Arr}, the following holds.
\begin{theorem}\label{thm:Tor>LArr}
Consider two categories ${\mathcal A}$ and ${\mathcal B}$ both of which
possess equalizers, and an object $({\mathcal T},(N_A,R_A),(N_B,R_B),\tau)$ in
$\mathrm{PreTor}({\mathcal A},{\mathcal B})$, that satisfies the assumptions in
Theorem \ref{thm:Tor>Arr}. Then the following assertions hold.

(1) The equalizer $ \left( D,j\right)
  =\mathrm{Equ}_{\mathrm{Fun}}(\theta^{l},\theta^{r})$ of the natural morphisms
$$
\theta^{l}:=(R_{B}\epsilon^A N_{B}R_{B}N_{A}R_{A} {N_{B}})
\circ ({R_{B}N_{A}}\tau)\qquad \textrm{and}\qquad
\theta ^{r}:=\eta^B R_{B}N_{A} R _{A}N_{B}
$$
defines a comonad ${\mathbb D}=(D,\Delta^D,\varepsilon^D)$ on ${\mathcal B}$
such that the functor $D$ preserves equalizers.

(2) There is an object
$({\mathcal T},(N_B,R_B),(N_A,R_A),{\mathbb D},\zeta)\in
  \overline{\mathrm{RArr}}({\mathcal
  B},{\mathcal A})$ where the comonad ${\mathbb D}$ was constructed in part (1).

(3) For a morphism $F$ between pre-torsors
  $({\mathcal T},(N_A,R_A),(N_B,R_B),\tau)$ and $({\mathcal
  T}',$ $(N'_A,R'_A),(N'_B,R'_B),\tau')$,
  both of which satisfy the conditions in Theorem \ref{thm:Tor>Arr}, there
  exists a unique morphism $t:{\mathbb D} \to {\mathbb D}'$ between
       the comonads in part (1) such that $(F,t)$ is a morphism
       between the objects of $\overline{\mathrm{RArr}}({\mathcal
  B},{\mathcal A})$ in part (2).
\end{theorem}

\subsection{The category equivalence}\label{sec:cat_eq}

In Section \ref{sec:Gal>Tor} we constructed functors from the categories
$\mathrm{RArr}({\mathcal A},{\mathcal B})$ and
$\overline{\mathrm{RArr}}({\mathcal
  B},{\mathcal A})$ to $\mathrm{PreTor}({\mathcal
    A},{\mathcal B})$. In Section \ref{sec:Tor>Gal} we constructed functors
  in the opposite direction, from an appropriate subcategory of pre-torsors to
  the categories $\mathrm{RArr}({\mathcal A},{\mathcal B})$ and
  $\overline{\mathrm{RArr}}({\mathcal B},{\mathcal A})$. The aim of this
  section is to find subcategories in
both categories, in such a way that the functors in Section \ref{sec:Gal>Tor}
and Section \ref{sec:Tor>Gal} establish equivalences between them.

Motivated by the premises of Theorem \ref{thm:Tor>Arr}, we impose the
following definition.

\begin{definition}\label{def:reg}
Let ${\mathcal A}$ and ${\mathcal B}$ be categories with equalizers and let
${\mathcal T}$ be any category.
An adjunction $(N:{\mathcal A} \to {\mathcal T}, R:{\mathcal T}\to {\mathcal
  A})$  is called {\em regular} whenever the following conditions hold.
\begin{itemize}
\item The unit $\eta$ of the adjunction is a regular natural monomorphism;
\item the functor $N$ preserves equalizers.
\end{itemize}
An object $({\mathcal T},(N_A,R_A),(N_B,R_B),\tau)$ of
$\mathrm{PreTor}({\mathcal A},{\mathcal B})$ is said to be {\em regular} if both
adjunctions $(N_A,R_A)$ and $(N_B,R_B)$ are regular.
We denote the full subcategory of regular pre-torsors by
$\mathrm{PreTor}_{reg}({\mathcal A}, {\mathcal B})$.

An object $({\mathcal T},(N_A,R_A),(N_B,R_B), {\mathbb C},\xi)$ of
$\mathrm{Arr}({\mathcal A}, {\mathcal B})$ is said to be
{\em  (co)-regular} if
$(R_A,\xi)$ is a (co)-regular comonad arrow,
both adjunctions $(N_A,R_A)$ and $(N_B,R_B)$ are regular and
the functor underlying ${\mathbb C}$ preserves equalizers.
We denote the full subcategories of regular and co-regular objects in
$\mathrm{Arr}({\mathcal A}, {\mathcal B})$
by $\mathrm{RArr}_{reg}({\mathcal A}, {\mathcal B})$ and
$\overline{\mathrm{RArr}}_{reg}({\mathcal A}, {\mathcal B})$, respectively.
\end{definition}

By corestriction, Theorem \ref{thm:Tor>Arr} provides us with a functor
$\Gamma: \mathrm{PreTor}_{reg}({\mathcal A}, {\mathcal B})\to
\mathrm{RArr}_{reg} ({\mathcal A}, {\mathcal B})$. By restriction and
corestriction, Theorem \ref{thm:Arr>Tor} yields a functor
$\Omega:\mathrm{RArr}_{reg} ({\mathcal A}, {\mathcal B}) \to
\mathrm{PreTor}_{reg} ({\mathcal A}, {\mathcal B})$.
Symmetrically, by Theorem \ref{thm:Tor>LArr} and Theorem \ref{thm:LArr>Tor},
we have functors between the categories $\mathrm{PreTor}_{reg}({\mathcal A},
{\mathcal B})$ and $\overline{\mathrm{RArr}}_{reg} ({\mathcal B}, {\mathcal
  A})$.
Our main result states that both pairs of functors are inverse equivalences.

\begin{theorem}\label{theo:II}
For two categories ${\mathcal A}$ and ${\mathcal B}$ with equalizers,
the following categories are equivalent.
\begin{itemize}
\item[{(i)}] $\mathrm{PreTor}_{reg}({\mathcal A}, {\mathcal B})$;
\item[{(ii)}] $\mathrm{RArr}_{reg} ({\mathcal A}, {\mathcal B})$;
\item[{(iii)}] $\overline{\mathrm{RArr}}_{reg} ({\mathcal B}, {\mathcal A})$.
\end{itemize}
\end{theorem}

\begin{proof}
{\em Equivalence of (i)and (ii)}.
Consider the functors $\Gamma$ and $\Omega$ in the paragraph preceding the
theorem.
First we construct a natural isomorphism $\mathrm{RArr}_{reg} ({\mathcal
  A}, {\mathcal B})\to \Gamma \Omega$. That is, we associate an isomorphism
\begin{eqnarray*}
({\mathcal T},N_A,N_B,w):&&({\mathcal T},(N_A,R_A),(N_B,R_B),{\mathbb C},
  \xi) \to \\
&&\Gamma \Omega
({\mathcal T},(N_A,R_A),(N_B,R_B),{\mathbb C}, \xi)=:
({\mathcal T},(N_A,R_A),(N_B,R_B),{\widetilde {\mathbb C}}, {\widetilde \xi})
\end{eqnarray*}
to any object
$({\mathcal T},(N_A,R_A),(N_B,R_B),{\mathbb C}, \xi)$
of $\mathrm{RArr}_{reg} ({\mathcal A}, {\mathcal B})$.
The functor ${\widetilde C}$, underlying the comonad ${\widetilde {\mathbb
C}}$, is defined as the equalizer of the morphisms $\omega^l$ and $\omega^r$ in
Theorem \ref{thm:Tor>Arr} (1), associated to the pre-torsor in Theorem
\ref{thm:Arr>Tor}. 
Equivalently, ${\widetilde C}$ is the equalizer of
\begin{eqnarray*}
&&(\xi^{-1} N_A R_A N_A) \circ \omega^l = (C \eta^A R_A N_A) \circ (\xi^{-1}
  N_A)\qquad \textrm{and}\\
&&(\xi^{-1} N_A R_A N_A) \circ \omega^r = (C R_A N_A \eta^A) \circ (\xi^{-1}
  N_A).
\end{eqnarray*}
On the other hand,
by regularity of $({\mathcal T},(N_A,R_A),(N_B,R_B),{\mathbb C}, \xi)$,
the functor $C$ underlying the comonad ${\mathbb C}$ preserves the equalizer
$\eta^A$. So we conclude that also
\begin{eqnarray*}
\left(C,  (\xi  N_A)\circ (C\eta^A)\right)&=&
\mathrm{Equ}_{\mathrm{Fun}}((C \eta^A R_A N_A) \circ (\xi^{-1} N_A),
(C R_A N_A \eta^A) \circ (\xi^{-1} N_A))\\
&=&\mathrm{Equ}_{\mathrm{Fun}}(\omega^l,\omega^r).
\end{eqnarray*}
Thus by uniqueness of an equalizer up to isomorphism,
there exists a natural isomorphism $w:C \to {\widetilde C}$, such
that the natural monomorphism ${\widetilde i}: {\widetilde C} \to
R_A N_B R_B N_A$ satisfies
\begin{equation}\label{eq:w}
{\widetilde i}\circ w= (\xi N_A) \circ (C\eta^A).
\end{equation}
The coproduct {$\Delta ^{\widetilde{C}}$} of the
comonad ${\widetilde {\mathbb C}}$ is defined as the unique
morphism satisfying \eqref{form:DeltaCbar}, i.e. the equation
$$
(C\eta^A C\eta^A)\circ (w^{-1} w^{-1}) \circ {\Delta
^{\widetilde{C}}} = (C\eta^A C\eta^A)\circ \Delta^C \circ w^{-1}.
$$
Since $C$ preserves equalizers by the regularity assumption,
$C\eta^A C\eta^A$ is a monomorphism. Thus we conclude that
$(ww)\circ\Delta^C ={\Delta ^{\widetilde{C}}} \circ w$.
Similarly, the counit ${\varepsilon ^{\widetilde{C}}}$
of ${\widetilde{\mathbb
    C}}$ is the unique morphism satisfying \eqref{form:epsilonC}, that is, the
equality
$$
\eta^A \circ {\varepsilon ^{\widetilde{C}}} = \eta^A
\circ \varepsilon^C \circ w^{-1}.
$$
Therefore, by monomorphism property of $\eta^A$, it follows that
$\varepsilon^C = {\varepsilon ^{\widetilde{C}}}\circ
w$, so $w$ is an isomorphism of comonads, as required.

By \eqref{eq:w}, naturality and an adjunction  relation, $\xi = {\widetilde
  \xi} \circ (w R_A)$. That is, $({\mathcal T},N_A,N_B,w)$ is a morphism in
$\mathrm{RArr}_{reg} ({\mathcal A}, {\mathcal B})$.

It remains to prove that the isomorphism $({\mathcal T},N_A,N_B,w)$ is
natural. That is,
given a morphism $(F,t):({\mathcal T},(N_A,R_A),(N_B,R_B),{\mathbb C},
\xi) \to ({\mathcal T}',(N'_A,R'_A),(N'_B,R'_B),{\mathbb C}',$ $ \xi')$
in $\mathrm{RArr}_{reg}({\mathcal A},{\mathcal B}) $ (with corresponding
natural morphisms $a$ and $b$ in \eqref{eq:a&b}),
the commutativity condition
\begin{equation}\label{eq:w_nat}
w' \circ t = {\widetilde t} \circ w
\end{equation}
holds, where ${\widetilde t}:\widetilde{C}\rightarrow
\widetilde{C^{\prime }}$ is defined by $(F,{\widetilde t}):=\Gamma
\Omega(F,t)$.
In order to prove \eqref{eq:w_nat}, compose both sides on the left by the
canonical monomorphism ${\widetilde i}':{\widetilde C}'\to R_A' N_B' R_B'
N_A'$. The resulting equivalent condition follows by the construction of
${\widetilde t}$ via the identity ${\widetilde i}'  \circ {\widetilde t} =
(R_A b R_B a) \circ {\widetilde i}$, \eqref{eq:w} and the fact
that $(F,t)$ is a morphism in $\mathrm{RArr}({\mathcal A}, {\mathcal B})$.
This completes the
proof of the claim that the functor $\Gamma \Omega$ is naturally isomorphic
to the identity functor $\mathrm{RArr}_{reg}({\mathcal A}, {\mathcal B})$.

In the converse order, we claim that $\Omega\Gamma$ is equal to the identity
functor on $\mathrm{PreTor}_{reg}({\mathcal A}, {\mathcal B})$. For an object
$({\mathcal T},(N_A,R_A),(N_B,R_B),\tau)$ in
$\mathrm{PreTor}_{reg}(\mathcal{A},\mathcal{B})$, denote
$$
({\mathcal T},(N_A,R_A),(N_B,R_B),{\widehat \tau}) :=
\Omega \Gamma ({\mathcal T},(N_A,R_A),(N_B,R_B),\tau).
$$
By Theorem \ref{thm:Arr>Tor} (1) and Theorem \ref{thm:Tor>Arr} (2),
$$
{\widehat \tau}=(\xi N_A R_A N_B)\circ (C\eta^A R_A N_B) \circ (\xi^{-1} N_B)
\circ (R_A N_B \eta^B),
$$
where $\xi$ is the isomorphism constructed in Lemma \ref{lem:xi}.
Using the explicit form of $\xi$ and the relation defining
its inverse, naturality and the adjunction
relations, ${\widehat \tau}$ is checked to be equal to ${\tau}$.
On the morphisms $\Omega \Gamma$ obviously acts as the identity
map.

Equivalence of (i) and (iii) is proven symmetrically.
\end{proof}

We can apply the equivalence functor obtained between
$\mathrm{PerTor}_{reg}({\mathcal A},{\mathcal B})$ and
$\mathrm{RArr}_{reg}({\mathcal A},{\mathcal B})$ to objects of the subcategory
$\mathrm{Gal}_{reg}({\mathcal A},{\mathcal B})$ of
$\mathrm{RArr}_{reg}({\mathcal A},{\mathcal B})$, cf. Section \ref{sec:entw}.

\begin{corollary}\rm
(1)
Consider a mixed distributive law of a monad ${\mathbb N}$ and a comonad
${\mathbb C}$ on a category ${\mathcal A}$ and denote the corresponding
lifting of ${\mathbb
C}$ to a comonad on ${\mathcal A}_{\mathbb N}$ by ${\widetilde {\mathbb
C}}$. Let $L$ be a ${\widetilde {\mathbb C}}$-Galois functor with right
adjoint $R$. Assume that both
adjunctions $(F_{\mathbb N},U_{\mathbb N})$ and $(L,R)$ are regular and the
functor underlying ${\mathbb C}$ preserves equalizers. Then there is a regular
pre-torsor
\begin{equation}\label{eq:entw.tor}
\big( {\mathcal A}_{\mathbb N}, (F_{\mathbb N},U_{\mathbb N}), (L,R),\tau\big).
\end{equation}
(2)
Conversely, consider a monad ${\mathbb N}$ on a category ${\mathcal A}$ and an
adjoint pair of functors $(L:{\mathcal B}\to {\mathcal A}_{\mathbb N},
  R:{\mathcal A}_{\mathbb N}\to {\mathcal B})$.
Then any regular pre-torsor of the form \eqref{eq:entw.tor} determines a
comonad ${\mathbb C}$ on ${\mathcal A}$ such that
\begin{itemize}
\item the functor underlying ${\mathbb C}$ preserves equalizers,
\item ${\mathbb C}$ lifts to a comonad ${\widetilde {\mathbb C}}$ on ${\mathcal
    A}_{\mathbb N}$,
\item $L$
carries the structure of a ${\widetilde {\mathbb C}}$-Galois functor.
\end{itemize}
\end{corollary}

More specifically, we can consider the examples in
Section \ref{sec:bimod.ex}.
Consider an entwining structure $({\bf T},{\bf
C},\boldsymbol{\psi})$ over an associative and unital algebra
${\bf A}$. Denote the induced ${\bf T}$-coring ${\bf T\ot_A C}$ by
${\widetilde {\bf C}}$. Then any Galois right ${\widetilde {\bf
C}}$-comodule $\Sigma$ determines a (not
necessarily regular) object in $\mathrm{RArr}(\mathrm{Mod}$-${\bf
A}, \mathrm{Mod}$-${\bf B})$, hence by Theorem \ref{thm:Arr>Tor} also a
pre-torsor
\begin{equation}\label{eq:bim.tor}
\big( \mathrm{Mod}\textrm{-}{\bf T}, ((-)\ot_{\bf A} {\bf
T},\mathrm{Hom}_{\bf T}({\bf T},-)), ((-)\ot_{\bf B}
\Sigma,\mathrm{Hom}_{\bf T}(\Sigma,-)), \tau\big),
\end{equation}
which is regular provided that ${\bf T}$ is a faithfully flat left
${\bf A}$-module and $\Sigma$ is a faithfully flat left ${\bf
  B}=\mathrm{End}^{\widetilde {\bf C}} (\Sigma)$-module.
Note, however, that not every (regular) pre-torsor over the given
adjunctions in \eqref{eq:bim.tor} arises from an entwining structure.

 For two algebras
${\bf A}$ and ${\bf B}$, consider an ${\bf
  A}$-ring ${\bf T}$ and a ${\bf B}$-${\bf T}$ bimodule $\Sigma$. Then any
regular pre-torsor of the form \eqref{eq:bim.tor} determines a comonad
${\mathbb C}$ on $\mathrm{Mod}\textrm{-}{\bf A}$ such that
\begin{itemize}
\item the functor underlying ${\mathbb C}$ preserves equalizers,
\item ${\mathbb C}$ lifts to a comonad ${\widetilde {\mathbb C}}$ on
  $\mathrm{Mod}\textrm{-}{\bf T}$,
\item $(-)\ot_{\bf B}\Sigma$ possesses
a ${\widetilde {\mathbb C}}$-Galois functor structure.
\end{itemize}
However, the comonad ${\mathbb C}$ is not known to be induced by a coring
(i.e. the underlying functor is not known to be a left adjoint).

Let us explain here the most important difference between the approaches to
the relation between Galois comodules and pre-torsors in the current paper on
one hand, and in \cite{BT} and \cite{BV} on the other hand.

In \cite[Theorem 3.4]{BT} pre-torsors of the form \eqref{eq:bim.tor} are
considered in the particular
case when the ${\bf B}$-${\bf T}$ bimodule $\Sigma$ is
a free rank 1 right ${\bf T}$-module (hence ${\bf T}$ is also a ${\bf
B}$-ring).
More generally, in \cite[Theorem 2.16]{BV} the case is discussed when
$\Sigma$ is a finitely generated and projective right ${\bf T}$-module.
Following the methods in \cite[Theorem 2.16]{BV}, in the case when $\Sigma$ is
in addition a faithfully flat right ${\bf A}$-module and ${\bf T}$ is a
faithfully flat left ${\bf A}$-module, a pre-torsor of
the form \eqref{eq:bim.tor} can be shown to determine an {\em ${\bf A}$-coring}
${\bf C}$.
However, note that the assumptions made on $\Sigma$ in this approach, do not
imply the regularity conditions in Definition \ref{def:reg}.

On the other hand, following Theorem \ref{thm:Tor>Arr} (1), one
can assume that in \eqref{eq:bim.tor} $\Sigma$ is a faithfully
flat left ${\bf B}$-module, and ${\bf T}$ is a faithfully flat
left ${\bf A}$-module. Under these assumptions we associated a
{\em comonad} $\mathbb C$ on $\mathrm{Mod}\textrm{-}{\bf A}$ to a
pre-torsor \eqref{eq:bim.tor}. Thus in the cases when both groups
of assumptions hold, both the coring ${\bf C}$ {\em and} the
comonad ${\mathbb C}$ can be constructed. However, there is no
reason to expect that the comonad ${\mathbb C}$ is induced by the
coring ${\bf C}$: The comonad ${\mathbb C}$ is unique
(up to a natural isomorphism) with the property that
the underlying functor preserves equalizers -- a property the
comonad $(-)\ot_{\bf A} {\bf C}$ is not known to obey (unless
${\bf C}$ is a flat left ${\bf A}$-module, e.g.
because we work with equal commutative base rings ${\bf A}$ and
${\bf B}$ and their symmetrical modules, cf. \cite[Remark 4.7]{BT}).

This deviation between the constructions in Section \ref{sec:pre-tor} on one
hand, and in the works  \cite{Hobst:phd}, \cite{BT} and \cite{BV} on the other
hand, shows that there is a conceptual ambiguity how to generalize faithfully
flat Hopf bi-Galois objects to non-commutative base algebras. Following the
(more conventional) approach in \cite{Hobst:phd}, \cite{BT} and \cite{BV}, the
coacting symmetry objects can be described by {\em two corings}. However, as
it was observed in \cite{BT}, Morita Takeuchi equivalence of these corings can
not be proven in general. Here we would like to point out an alternative
strategy: One can allow for the coacting symmetry structures to be {\em two
comonads}, whose underlying functors are not necessarily left adjoints but, as
a gain, they preserve kernels. As it is proven in Section \ref{sec:Mor-Tak},
in this setting Morita Takeuchi equivalence of the two comonads is easily
proven.

\section{Equivalence of comodule categories}
\label{sec:Mor-Tak}

In Section \ref{sec:pre-tor} we proved equivalences between three categories
$\mathrm{RArr}_{reg}({\mathcal A},{\mathcal B})$, \break
$\overline{\mathrm{RArr}}_{reg}({\mathcal B}, {\mathcal A})$ and
$\mathrm{PreTor}_{reg} ({\mathcal A}, {\mathcal B})$, for two categories
${\mathcal A}$ and ${\mathcal B}$ possessing equalizers. In this way we
associated in particular two comonads, ${\mathbb C}$ on ${\mathcal A}$ and
${\mathbb D}$ on ${\mathcal B}$, to any object of $\mathrm{PreTor}_{reg}
({\mathcal A}, {\mathcal B})$.
In this section we prove that the comonads ${\mathbb C}$ and
${\mathbb D}$ have equivalent comodule categories, what generalizes the result
in \cite{Sch} about Morita Takeuchi equivalence of the two Hopf algebras
associated to a bi-Galois object.
The proof is presented at the level of generality in Lemma \ref{lem:C&D}.

Note that, for any comonad ${\mathbb D}$ on a category ${\mathcal B}$, the
forgetful functor $U^{\mathbb D}$ is a left ${\mathbb D}$-comodule functor via
the coaction $U^{\mathbb D}\gamma^{D}$
(notation introduced in Definition \ref{def:monad} (6)). In Proposition
\ref{prop:bicomod>functor}, a {\em `cotensor product'} of $U^{\mathbb D}$ with 
a right ${\mathbb D}$-comodule functor occurs.

\begin{proposition}\label{prop:bicomod>functor}
Let ${\mathcal A}$ and ${\mathcal B}$ be categories possessing equalizers,
${\mathbb C}=(C,\Delta^C,\varepsilon^C)$ be a comonad on ${\mathcal A}$,
${\mathbb D}=(D,\Delta^D,\varepsilon^D)$ be a comonad on
${\mathcal B}$ and $(Q,c,d)$ be a ${\mathbb C}$-${\mathbb D}$ bicomodule
functor. If the functor $C$ preserves equalizers, then
there is a functor $I_Q:{\mathcal B}^{\mathbb D} \to {\mathcal A}^{\mathbb C}$
such that
$$
\xymatrix{
U^{\mathbb C} I_Q \ar[r]&
Q U^{\mathbb D} \ar@<2pt>[rr]^-{d U^{\mathbb D}}\ar@<-2pt>[rr]_-{Q U^{\mathbb D}
  \gamma^D} &&
QDU^{\mathbb D}
}
$$
is an equalizer.
\end{proposition}

\begin{proof}
Since ${\mathcal A}$ is assumed to have equalizers, by Lemma \ref{Lem: Equ}
there exists the equalizer
$$
\xymatrix{
I_0 \ar[r]^-{e}&
Q U^{\mathbb D} \ar@<2pt>[rr]^-{d U^{\mathbb D}}\ar@<-2pt>[rr]_-{Q U^{\mathbb D}
  \gamma^D} &&
QDU^{\mathbb D}
}.
$$
We construct $I_Q:{\mathcal B}^{\mathbb D} \to {\mathcal A}^{\mathbb C}$ by
equipping $I_0$ with a left ${\mathbb C}$-coaction $c_0:I_0 \to C I_0$ and
putting $I_Q X:= \big( I_0 X, c_0 X)$, for any object $X$ in ${\mathcal
B}^{\mathbb D}$, and $I_Q f := I_0 f$, for any morphism $f$.

By naturality, the definition of $e$ and the bicomodule property of $Q$, we
conclude that 
$$
(CQU^{\mathbb D} \gamma^D) \circ (c U^{\mathbb D}) \circ e =
(C d U^{\mathbb D}) \circ (c U^{\mathbb D}) \circ e.
$$
Since $Ce$ is the equalizer of $C d U^{\mathbb D}$ and $CQ U^{\mathbb D}
  \gamma^D$, its universality implies the existence of a natural
morphism $c_0:I_0 \to C I_0$ such that
\begin{equation}\label{eq:c.0}
(Ce) \circ c_0 = (c U^{\mathbb D})\circ e.
\end{equation}
Its coassociativity and counitality are immediate by
coassociativity and counitality of $c$.
\end{proof}

In light of Proposition \ref{prop:bicomod>functor}, a functor $Q$ in Lemma
\ref{lem:C&D} induces a functor $I_Q:{\mathcal B}^{\mathbb D} \to {\mathcal
  A}^{\mathbb C}$ for the comonads ${\mathbb C}$ and ${\mathbb D}$ constructed
in parts (1) and (2) of Lemma \ref{lem:C&D}, respectively. Our next task is to
construct a  ${\mathbb D}$-${\mathbb C}$ bicomodule functor
${\overline Q}$ such that
the induced functor $I_{\overline Q}:{\mathcal A}^{\mathbb C} \to {\mathcal
  B}^{\mathbb D}$ yields the inverse of $I_Q$.

\begin{lemma}\label{lem:Q.bar}
In the setting of Lemma \ref{lem:C&D} and using the notations
in its proof, define functors ${\overline Q}$ and $\overline{Q}':{\mathcal A}
\to {\mathcal B}$ via the respective equalizers
$$
\xymatrix{
{\overline Q}\ar[r]^-{q}&
PC \ar@<2pt>[rr]^-{(\theta^l P)\circ (Pi)}
\ar@<-2pt>[rr]_-{(\theta^r P)\circ (Pi)}&&
SPQP
} \quad \textrm{and}\quad
\xymatrix{
{\overline Q}'\ar[r]^-{q'}&
DP \ar@<2pt>[rr]^-{(P\omega^l)\circ (jP)}
\ar@<-2pt>[rr]_-{(P \omega^r)\circ (jP)}&&
PQPR}.
$$
The following statements hold.

(1) The functors ${\overline Q}$ and ${\overline Q}'$ are naturally
isomorphic (hence they can be chosen equal).

(2) The functor ${\overline Q}$ can be equipped with the structure of a
${\mathbb D}$-${\mathbb C}$ bicomodule.
\end{lemma}

\begin{proof}
(1)
By construction, $(jP)\circ q'$ equalizes $P\omega^l$ and $P\omega^r$.
Thus by universality of the equalizer $Pi$, there exists a natural morphism
$\nu_0:\overline{Q}'\to PC$ such that
\begin{equation}\label{eq:nu_0}
(Pi)\circ \nu_0 = (jP)\circ q'.
\end{equation}
Moreover, since $j$ equalizes $\theta^l$ and $\theta^r$ by definition,
$(jP)\circ q'$ equalizes $\theta^l P$ and $\theta^r P$.
Therefore
\eqref{eq:nu_0} implies that $\nu_0$ equalizes $(\theta^l P)\circ (Pi)$ and
$(\theta^r P)\circ (Pi)$. Thus by universality of the equalizer $q$, there is
a natural morphism $\nu: \overline{Q}'\to \overline{Q}$ such that $q \circ \nu
= \nu_0$. Equivalently,
\begin{equation}\label{eq:nu}
(Pi)\circ q \circ \nu = (jP)\circ q'.
\end{equation}
By a symmetrical reasoning, there is a morphism $\nu': \overline{Q}\to
\overline{Q}'$, such that
\begin{equation}\label{eq:nu'}
(jP)\circ q' \circ \nu' = (Pi)\circ q.
\end{equation}
By \eqref{eq:nu} and \eqref{eq:nu'}, $(jP)\circ q'\circ \nu' \circ \nu=
(jP)\circ q'$. Since $jP$ and $q'$ are monomorphisms, this implies $\nu' \circ
\nu = \overline{Q}'$. A symmetrical reasoning justifies $\nu \circ \nu'
=\overline{Q}$.

Since an equalizer is defined up to isomorphism, we may choose ${\overline
  Q}={\overline Q}'$, resulting in the identity
\begin{equation}\label{eq:q.q'}
(jP)\circ q' = (Pi)\circ q .
\end{equation}

(2) By assumption (ii) in Lemma \ref{lem:C&D} and naturality,
$$
(\theta^l PQP)\circ (P \tau P) = (SP\tau P)\circ (\theta^l P)
\qquad \textrm{and} \qquad
(\theta^r PQP)\circ (P \tau P) = (SP\tau P)\circ (\theta^r P) .
$$
Together with \eqref{eq:q.q'} and \eqref{form:DeltaCbar}, this implies
\begin{equation}\label{eq:st.2}
(SPQPi) \circ (\theta^l PC) \circ (PiC) \circ (P\Delta^C ) \circ q =
(SPQPi) \circ (\theta^r PC) \circ (PiC) \circ (P\Delta^C ) \circ q .
\end{equation}
Since $SPQPi$ is a monomorphism, we conclude by universality of the equalizer
$qC$ that there exists a unique morphism ${\overline c}:{\overline Q}\to
{\overline  Q} C$ such that
\begin{equation}\label{eq:c.bar}
(qC) \circ {\overline c} = (P \Delta^C) \circ q.
\end{equation}
Obviously, ${\overline c}$ is a coassociative and counital coaction.
Symmetrically, a ${\mathbb D}$-coaction $\overline{d}:\overline{Q} \to D
\overline{Q}$ is defined by the condition
\begin{equation}\label{eq:d.bar}
(Dq')\circ \overline{d} = (\Delta^D P) \circ q' .
\end{equation}
The ${\mathbb C}$-, and $\mathbb D$-coactions on $\overline{Q}$ commute by
\eqref{form:DeltaCbar} and the analogous formula for $\Delta^D$ in
\eqref{eq:D.str}, assumption (ii) in Lemma \ref{lem:C&D} and \eqref{eq:q.q'}.
\end{proof}

\begin{theorem}\label{thm:comod.eq}
In the setting of Lemma \ref{lem:C&D}, consider the ${\mathbb C}$-${\mathbb D}$
bicomodule $Q$ in  Lemma \ref{lem:C&D} (3) and the ${\mathbb D}$-${\mathbb C}$
bicomodule ${\overline Q}$ in Lemma \ref{lem:Q.bar}. Then the induced functors
$I_Q:{\mathcal B}^{\mathbb D} \to {\mathcal A}^{\mathbb C}$ and
$I_{\overline Q}:{\mathcal A}^{\mathbb C} \to {\mathcal B}^{\mathbb D}$
(cf. Proposition \ref{prop:bicomod>functor}) are inverse equivalences.
\end{theorem}

\begin{proof}
Recall that
$I_{\overline 0}= U^{\mathbb D} I_{\overline Q}$ fits the equalizer
\begin{equation}\label{eq:I.0.bar}
\xymatrix{
I_{\overline 0} \ar[r]^-{\overline e}&
{\overline Q} U^{\mathbb C} \ar@<2pt>[rr]^-{{\overline c} U^{\mathbb C}}
\ar@<-2pt>[rr]_-{{\overline Q} U^{\mathbb C} \gamma^C} &&
{\overline Q} C U^{\mathbb C},
}
\end{equation}
where the ${\mathbb C}$-coaction ${\overline c}$ on ${\overline Q}$ was
constructed in \eqref{eq:c.bar}.

First we construct a natural isomorphism between $I_{\overline 0} I_Q$ and
$U^{\mathbb D}$ and then show that it lifts to a natural isomorphism between
$I_{\overline Q} I_Q$ and the identity functor ${\mathcal B}^{\mathbb D}$.
By \eqref{eq:I.0.bar},
$$
({\overline c} I_0)\circ ({\overline e} I_Q) = ({\overline Q} c_0) \circ
({\overline e} I_Q),
$$
where the ${\mathbb C}$-coaction $c_0 = U^{\mathbb C} \gamma^C I_Q$ on $I_0$
was constructed in
\eqref{eq:c.0}. Compose both sides of this equality on the left by
$(P\varepsilon^C ie) \circ (qC I_0)$, and use \eqref{eq:c.bar}, \eqref{eq:c.0}
and \eqref{eq:C-coac} to conclude that
\begin{equation}\label{eq:auxi}
(PiQU^{\mathbb D}) \circ (qe) \circ ({\overline e} I_Q) =
(P \tau U^{\mathbb D}) \circ (P \varepsilon^C Q U^{\mathbb D}) \circ  (qe)
\circ ({\overline e} I_Q).
\end{equation}
By the counit formula in \eqref{eq:D.str}, \eqref{eq:q.q'}, fork property of
\eqref{eq:C.eq}, naturality, \eqref{eq:q.q'} again, fork property of
\eqref{eq:D.eq}, \eqref{eq:q.q'} and \eqref{form:epsilonC}, it follows that
$(sPr)\circ (\varepsilon^D P)\circ q'=(sPr)\circ (P \varepsilon^C)\circ
q$. Since $sPr$ is a monomorphism, this implies
$$
(\varepsilon^D P)\circ q'= (P \varepsilon^C)\circ q.
$$
These identities, together with \eqref{eq:q.q'}, \eqref{form:epsilonC} and the
analogous formula for $\varepsilon^D$ in \eqref{eq:D.str} imply that
$\alpha_0:= (P\varepsilon^C Q U^{\mathbb D}) \circ (qe) \circ ({\overline e}
I_Q)$ equalizes $\theta^l U^{\mathbb D}$ and $\theta^r U^{\mathbb D}$. Thus by
universality of the
equalizer $j U^{\mathbb D}$ (cf. Lemma \ref{Lem: Funct1}),
there exists a natural morphism $\alpha:
I_{\overline 0} I_Q\to D U^{\mathbb D}$ such that $(j U^{\mathbb D}) \circ
\alpha = \alpha_0$.
Note that by \eqref{eq:auxi}, $(P\tau U^{\mathbb D})\circ \alpha_0$ is a
monomorphism. This implies that $\alpha_0$ is a monomorphism and hence
$\alpha$ is a monomorphism.

By the coproduct formula in \eqref{eq:D.str}, \eqref{eq:D.coac} and
the definition of $e$, it follows that
$$
(jj U^{\mathbb D})\circ (\Delta^D U^{\mathbb D}) \circ \alpha =
(jj U^{\mathbb D})\circ (DU^{\mathbb D} \gamma^D) \circ \alpha.
$$
Since $jj U^{\mathbb D}$ is a monomorphism and $U^{\mathbb D} \gamma^D$ is the
equalizer of $\Delta^D U^{\mathbb D}$ and $DU^{\mathbb D} \gamma^D$, there
exists a natural morphism $\widetilde{\beta}: I_{\overline 0} I_Q \to
U^{\mathbb D}$, such that
\begin{equation}\label{eq:beta.tilde}
(U^{\mathbb D} \gamma^D)\circ \widetilde{\beta} = \alpha.
\end{equation}
We prove that $\widetilde{\beta}$ is an isomorphism by constructing its
inverse. Consider the morphism $\beta_1:= (Pc U^{\mathbb D})\circ (j
U^{\mathbb D}) \circ (U^{\mathbb D} \gamma^D): U^{\mathbb D} \to PCQU^{\mathbb
  D}$. We claim that there exists a morphism $\beta:  U^{\mathbb D} \to
I_{\overline 0} I_Q $, such that $(qe) \circ ({\overline e} I_Q) \circ \beta =
\beta_1$. Equivalently,
\begin{equation}\label{eq:beta}
\alpha \circ \beta = U^{\mathbb D} \gamma^D.
\end{equation}
This will be verified in three steps. First observe that by \eqref{eq:D.coac},
\eqref{eq:C-coac}, assumption (ii) in Lemma \ref{lem:C&D} and coassociativity
of the ${\mathbb D}$-coaction $\gamma^D:U^{\mathbb D} \to D U^{\mathbb D}$, it
follows that
$$
(PiQjU^{\mathbb D}) \circ (PC d U^{\mathbb D})\circ \beta_1 =
(PiQjU^{\mathbb D}) \circ (PC Q U^{\mathbb D}\gamma^D)\circ \beta_1.
$$
Since $PiQjU^{\mathbb D}$ is monic, it follows by universality of the equalizer
$PCe$ that there exists $\beta_2: U^{\mathbb D} \to PC I_0$ such that $(P
C e) \circ \beta_2 = \beta_1$.
Next one checks with similar steps that $\beta_2$ equalizes $(\theta^l P
I_0)\circ (PiI_0)$ and $(\theta^r P I_0)\circ (PiI_0)$, hence by universality
of the equalizer $q I_0$ there exists $\beta_3: U^{\mathbb D} \to \overline{Q}
I_0$, such that $(qI_0)\circ \beta_3 = \beta_2$. Finally, $\beta_3$ is checked
to equalize $\overline{Q} c_0$ and $\overline{c} I_0$, hence by universality
of the equalizer $\overline{e} I_Q$ there exists $\beta:U^{\mathbb D} \to
I_{\overline 0} I_Q $ such that $(\overline{e} I_Q) \circ \beta =
\beta_3$. This morphism $\beta$ clearly satisfies \eqref{eq:beta}.

Since both $\alpha$ and $U^{\mathbb D} \gamma^D$ are monomorphisms,
\eqref{eq:beta.tilde} and \eqref{eq:beta} imply that $\beta$ is the inverse of
$\widetilde{\beta}$.

Note that composing both sides of it on the left by the
monomorphism $PiQU^{\mathbb D}$, \eqref{eq:beta} can be written equivalently
in the form
\begin{equation}\label{eq:beta.var}
(PiQ U^{\mathbb D}) \circ (qe) \circ (\overline{e} I_Q) \circ \beta =
(P \tau U^{\mathbb D}) \circ (j U^{\mathbb D}) \circ (U^{\mathbb D} \gamma^D).
\end{equation}
The natural isomorphism $I_{\overline Q} I_Q \cong {\mathcal B}^{\mathbb D}$
is proven by showing that $\beta$ is a ${\mathbb D}$-comodule morphism, i.e.
\begin{equation}\label{eq:beta.colin}
({\overline d}_0 I_Q)\circ \beta = (D \beta) \circ (U^{\mathbb D} \gamma^D),
\end{equation}
where the ${\mathbb D}$-coaction ${\overline d}_0 :I_{\overline 0} \to D
I_{\overline 0}$ is defined via the condition
\begin{equation}\label{eq:d.bar.0}
(D {\overline e}) \circ {\overline d}_0 = ({\overline d} U^{\mathbb C}) \circ
  {\overline e}.
\end{equation}
In order to prove \eqref{eq:beta.colin}, compose both sides of it on the left
by the
monomorphism $(jPiQU^{\mathbb D})\circ (Dqe) \circ (D {\overline e} I_Q)$, and
use \eqref{eq:d.bar.0}, \eqref{eq:q.q'}, \eqref{eq:d.bar}, \eqref{eq:D.str},
\eqref{eq:beta.var}, assumption (ii) in Lemma \ref{lem:C&D} and then again
\eqref{eq:D.str} and \eqref{eq:beta.var}.

The natural isomorphism $I_Q I_{\overline Q}\cong {\mathcal A}^{\mathbb C}$ is
verified by similar steps: First an isomorphism $I_0 I_{\overline Q} \cong
U^{\mathbb C}$ is constructed, and it is proven to lift to an isomorphism
$I_Q I_{\overline Q}\cong {\mathcal A}^{\mathbb C}$.
Consider the natural monomorphism $\nu_0:= (QP \varepsilon^C U^{\mathbb C})
\circ (Q q U^{\mathbb C}) \circ (Q \overline{e}) \circ (e I_{\overline Q}):I_0
I_{\overline Q} \to QPU^{\mathbb C}$. It is checked to equalize $\omega^l
U^{\mathbb C}$ and  $\omega^r U^{\mathbb C}$, hence by universality of the
equalizer $i U^{\mathbb C}$ it determines a monomorphism $\nu: I_0 I_{\overline
  Q} \to C U^{\mathbb C}$ such that $(i U^{\mathbb C})\circ \nu
=\nu_0$. Furthermore, $\nu$ is checked to equalize $\Delta^C U^{\mathbb C}$
and $C U^{\mathbb C} \gamma^C$. Since $U^{\mathbb C} \gamma^C$ is the
equalizer of $\Delta^C U^{\mathbb C}$ and $C U^{\mathbb C} \gamma^C$, there
exists a natural morphism ${\widetilde \kappa}: I_0 I_{\overline Q} \to
U^{\mathbb C}$ such that
\begin{equation}\label{eq:kappa.tilde}
(U^{\mathbb C} \gamma^C)\circ {\widetilde \kappa} =\nu.
\end{equation}
The inverse of ${\widetilde \kappa}$ is constructed in three steps. First the
morphism $\kappa_1:= (iC U^{\mathbb C})\circ (\Delta^C U^{\mathbb C}) \circ
(U^{\mathbb C} \gamma^C)$ is checked to equalize $(Q \theta^l P U^{\mathbb
  C})\circ (QP i U^{\mathbb C})$ and $(Q \theta^r P U^{\mathbb C})\circ (QP i
U^{\mathbb C})$, hence by universality of the equalizer $Qq U^{\mathbb C}$, it
determines a morphism $\kappa_2:U^{\mathbb C}\to Q \overline{Q} U^{\mathbb
  C}$, such that $(Qq U^{\mathbb C})\circ \kappa_2 = \kappa_1$. Next
$\kappa_2$ is checked to equalize $Q \overline{c}$ and $Q \overline{Q}
U^{\mathbb C} \gamma^C$. Hence by universality of the equalizer $Q
\overline{e}$, it determines a morphism $\kappa_3:U^{\mathbb C} \to Q
I_{\overline 0}$, such that $(Q \overline{e})\circ \kappa_3 =
\kappa_2$. Finally $\kappa_3$ is shown to equalize $d I_{\overline 0}$ and $Q
{\overline d}_0$ (where ${\overline d}_0:I_{\overline 0}\to D I_{\overline 0}$ is
the ${\mathbb D}$-coaction). Hence by universality of the equalizer $e
I_{\overline Q}$, it determines a morphism $\kappa:U^{\mathbb C}\to I_0
I_{\overline Q}$ such that $(e I_{\overline Q})\circ \kappa =
\kappa_3$. Equivalently,
\begin{equation}\label{eq:kappa}
\nu \circ \kappa = U^{\mathbb C} \gamma^C.
\end{equation}
We conclude by \eqref{eq:kappa.tilde} and \eqref{eq:kappa} that $\kappa$ and
$\widetilde{\kappa}$ are mutual inverses.
By very similar steps to those used in the case of $\beta$, also $\kappa$ is
checked to be a ${\mathbb C}$-comodule morphism, i.e. to satisfy $(c_0
I_{\overline Q})\circ \kappa = (C \kappa)\circ (U^{\mathbb C} \gamma^C)$. This
proves that $\kappa$ lifts to the stated isomorphism $I_Q I_{\overline Q}\cong
{\mathcal A}^{\mathbb C}$.
\end{proof}

The following corollary is immediate by Theorem \ref{thm:comod.eq}.

\begin{corollary}\label{cor:MorTak}
Consider an object in $\mathrm{PreTor}_{reg}({\mathcal
A},{\mathcal B})$ and its images
in $\mathrm{RArr}_{reg}({\mathcal A},{\mathcal B})$ and
$\overline{\mathrm{RArr}}_{reg}({\mathcal B},{\mathcal A})$,
respectively, under the equivalences in Theorem \ref{theo:II}.
Then the occurring comonads ${\mathbb C}$ on ${\mathcal A}$ and
${\mathbb D}$ on ${\mathcal B}$ have equivalent categories of comodules.
\end{corollary}

The ${\bf A}$-${\bf B}$ pre-torsors in \cite[Corollary 4.8]{BT} induce regular
pre-torsors in the sense of Definition \ref{def:reg}. Moreover, the associated
comonads on $\mathrm{Mod}\textrm{-}{\bf A}$ and $\mathrm{Mod}\textrm{-}{\bf
B}$ are induced by the corings in \cite[Theorem 3.4]{BT}. Therefore,
Corollary \ref{cor:MorTak} extends \cite[Corollary 4.8]{BT}.

\section*{Acknowledgements}
We would like to express our deepest gratitude to Viola Bruni for
her careful reading of preliminary versions of our paper and many
useful observations on them. We are indebted to Giulia Travasoni
for discussions at an early stage or our work. We are grateful for
their interest and inspiring comments to the participants of the
conference {\em  "Categorical methods for Rings and Modules''}
December 3-6, 2007, Murcia, Spain, where our results were first
presented. We would like to thank Gene Abrams for his
helpful advices. The first author acknowledges a Bolyai J\'anos
Research Scholarship and financial support of the Hungarian
Scientific Research Fund OTKA F67910 and of the
European Commission grant MKTD-CT-2004-509794. This paper was
written while the second author was a member of G.N.S.A.G.A. with
partial financial support from M.I.U.R.. She would like to express
her gratitude to the members of the Research Institute for
Particle and Nuclear Physics, Budapest, for supporting her and
also for a very warm hospitality during her visit to the Institute
when this research began.


\begin{thebibliography}{BBW}
\bibitem[ABM]{ABM} A. Ardizzoni, G. B\"ohm and C. Menini, {\em A
Schneider type theorem for Hopf algebroids}, J. Algebra {318} (2007)
225--269.
{\em Corrigendum}, to be published.
A corrected version is available also at {\tt arXiv:math/0612633v2}.

\bibitem[BW]{TTT} M. Barr and T. Wells, \emph{Toposes, Triples and
  Theories}, Grundl. der math. Wiss. 278, Springer-Verlag,
  1983, available at \\
  {\tt http://www.case.edu/artsci/math/wells/pub/ttt.html}.

\bibitem[Be]{Be} J. Beck, {\em Distributive laws}, in {Seminar on Triples
 and Categorical Homology Theory},
 B. Eckmann (ed.), Springer LNM 80, pp 119-140, 1969.

\bibitem[BB]{BT} G. Böhm and T. Brzezi\'{n}ski, \emph{Pre-torsors and
equivalences}, J. Algebra  317 (2007) 544-580. \emph{Corrigendum}, J. Algebra
  319 (2008) 1339--1340.

\bibitem[BrMa]{BrMa} T. Brzezi\'nski and S. Majid, {\em Coalgebra Bundles},
  Commun. Math. Phys. 191 (1998) 467-492.

\bibitem[BV]{BV} T. Brzezi\'nski and J. Vercruysse, {\em Bimodule herds},
  preprint arXiv:0805.2510v2.

\bibitem[MW]{MW} B.Mesablishvili and R. Wisbauer, {\em Bimonads and Hopf
  monads on categories}, preprint arXiv:0804.1460v3.

\bibitem[D]{Dubuc} E. Dubuc, \emph{Kan extensions in enriched category
  theory}. Lecture Notes in Mathematics 145, Springer Verlag Berlin 1970.

\bibitem[GT]{GT} J. Gómez-Torrecillas, \emph{Comonads and Galois corings},
Appl. Categorical Structures 14 (2006) 579-598.

\bibitem[G]{Gr} C. Grunspan, \emph{Quantum torsors}, J. Pure Appl. Algebra 184
  (2003) 229--255.

\bibitem[H]{Hobst:phd} D. Hobst, {\em Antipodes in the theory of noncommutative
  torsors}. PhD thesis Ludwig-Maximilians Universit\"at M\"unchen 2004, Logos
  Verlag Berlin, 2004.

\bibitem[J]{Johnst} P. T. Johnstone, \emph{Adjoint lifting theorems for
  categories of algebras}, Bull. London Math. Soc. 7 (1975) 294-297.

\bibitem[KS]{KeSt} G. M. Kelly  and R. Street,
  {\em Review of the elements of 2-categories},
  Category Sem., Proc., Sydney 1972/1973, Lect. Notes Math. 420 (1974)
  75-103.

\bibitem[S]{S} R. Street, {\em The formal theory of monads}, J. Pure
  Appl. Algebra 2 (1972) 149-168.

\bibitem[Sch1]{Scha:Qtor} P. Schauenburg, {\em Quantum torsors with fewer
  axioms}, preprint arXiv math.QA/0302003.

\bibitem[Sch2]{Sch} P. Schauenburg, \emph{Hopf bi-Galois extensions},
  Comm. Algebra 24 (1996) 3797--3825.

\bibitem[Sch3]{Scha:ba.nc.base} P. Schauenburg, {\em Bialgebras over
  noncommutative rings and a structure theorem for Hopf bimodules},
  Appl. Categ. Structures {6} (1998) 193--222.

\bibitem[Sch4]{Scha:HbiGal} P. Schauenburg, {\em Hopf-Galois and bi-Galois
  extensions}, in: Galois theory, Hopf algebras, and semiabelian categories,
  pp 469--515, Fields Inst. Commun., {43}, Amer. Math. Soc., Providence,
  RI, 2004.

\bibitem[\v Sk]{Skoda} Z. \v Skoda, {\em Quantum heaps, cops and heapy
  categories,} Mathematical Communications 12, (2007) 1-9.

\bibitem[W]{W} R. Wisbauer, \emph{Algebras versus coalgebras},
  Appl. Cat. Str. 16 (2008) 255-295.

\end{thebibliography}
\end{document}